\documentclass[12pt,reqno]{article}

\usepackage[usenames]{color}
\usepackage{amssymb}
\usepackage{graphicx}
\usepackage{amscd}

\usepackage[colorlinks=true,
linkcolor=webgreen,
filecolor=webbrown,
citecolor=webgreen]{hyperref}

\definecolor{webgreen}{rgb}{0,.5,0}
\definecolor{webbrown}{rgb}{.6,0,0}

\usepackage{color}
\usepackage{float}

\usepackage{graphics,amsmath,amssymb}
\usepackage{amsfonts}
\usepackage{latexsym}
\usepackage{epsf}

\begin{document}

\begin{center}
\vskip 1cm
{\LARGE\bf Enumeration of Certain Classes of\\[3pt] $T_0$-hypergraphs}\\
\bigskip
\Large
Goran Kilibarda\footnote{Department of Mathematics, ALFA University, Belgrade, Serbia
(\href{mailto:gkilibar@gmail.com}{\tt gkilibar@gmail.com})}
and
Vladeta Jovovi\'{c}\footnote{Faculty of Technology and Metallurgy, University of Belgrade, Karnegijeva 4, Belgrade, Serbia\\
(\href{mailto:vladeta@eunet.rs}{\tt vladeta@eunet.rs})}
\\
\end{center}

\thispagestyle{plain}
\pagestyle{myheadings}

\newtheorem{theorem}{Theorem}[section]
\newtheorem{proposition}{Proposition}[section]
\newtheorem{corollary}{Corollary}[section]
\newtheorem{lemma}{Lemma}[section]


\makeatletter

\vskip .3in

\begin{abstract}
  A hypergraph is a $T_0$-hypergraph if for every two different vertices of the hypergraph there exists an
  edge containing one of the vertices and not containing the other. A general method for the enumeration of
  certain classes of $T_0$-hypergraphs is given. $T_0$-hypergraphs that are considered here are singled out
  both by the properties they themselves satisfy and by the properties that dual hypergraphs associated with
  them satisfy. Though in case of the so-called ordered hipergraphs the property `to be a $T_0$-hypergraph'
  is reduced to the property `to having different columns' of corresponding matrices, combining this property
  with some properties, that we are considering here, gives sometimes classes of hypergraphs that are not so
  easy to enumerate. The problem of enumerating some of thus obtained classes remains unsolved. Special attention
  is devoted to enumerating of different classes of covers and connected hypergraphs.
\end{abstract}

\section{Introduction}

A non-empty finite set together with a finite family of its subsets is called a hypergraph. The elements of
the set are called vertices, and the members of the family are called edges of the given hypergraph. A hypergraph
could be labelled or unlabelled, and all hypergraphs in the paper are labelled. However, all the problems
investigated here could also be considered for the unlabelled case, though most of them are unsolved by now.
Just like in general topology we speak about $T_0$-spaces, here we speak about $T_0$-hypergraphs. A hypergraph
is a $T_0$-hypergraph if for every two different vertices there exists an edge which contains exactly one of
these vertices.

A hypergraph is ordered, freely speaking, if a linear order is given on the family of its edges, otherwise it is
unordered. To every ordered labelled hypergraph $H$ one can assign a binary matrix --- its incidence matrix $M_H$.
It is easy to see that an ordered labelled hypergraph $H$ is a $T_0$-hypergraph iff the matrix $M_H$ has no equal
columns. The dual hypergraph $H^T$ of an ordered labelled hypergraph $H$ is the hypergraph whose incidence matrix
is the transpose of $M_H$. Let a hypergraph property $\mathfrak{p}$ be given. We say that $H$ has the dual property
of $\mathfrak{p}$, that is, $H$ is a dual $\mathfrak{p}$-hypergraph, if $H^T$ has the property $\mathfrak{p}$. In
the paper we suppose that, in general, two sets $\mathfrak{P}$ and $\mathfrak{P}'$ of properties are given, and we
consider (both in an ordered and in an unordered case) the class of all hypergraphs having all the properties from
$\mathfrak{P}$ and having all the dual properties of the properties from $\mathfrak{P}'$. The properties we consider
are: `to be a cover', `to be $k$-dimensional', `to be $k$-uniform', `to be without intersecting property', etc., and,
clearly, all possible combinations of these properties. Our aim is to find how many hypergraphs of the class defined
in such a way have $T_0$-property. We suggest a more general method for enumerating such classes. Theorem \ref{trm1},
and its consequences Theorems \ref{trm2}, \ref{trm3}, \ref{trm4} and \ref{trm5}, constitute the base of the method.
For a larger part of the enumerations it is enough to apply Theorem \ref{trm4}, but there are some `hard nuts' for
which we are forced to use Theorem \ref{trm2}. In our opinion, Theorem \ref{42} and Theorem \ref{52} concerning covers
and connected hypergraph are also among the more significant theorems. Some of the problems adduced here remain unsolved
by now.

The results, obtained in this paper, concerning covers are related to corresponding results for $k$-covers that
are adduced and cited in \cite{Gou}. The results, concerning connected hypergraphs, are related to corresponding
results that are cited in \cite{Kre}.

The paper is completely self-contained, and all necessary notions and notations are given in Section 2. All the
formulas given here are tested by computer, and the majority of corresponding sequences, in an unordered case,
are published in \cite{Slo1}.

\section{Basic notions and notations}

Let $X$ be a set. Denote by
$|X|$
the cardinality of $X$. If
$|X|=n$,
then we say that $X$ is an $n$-set.

Let
$m_1,m_2\in\text{\bf Z}$, $m_1\le m_2$,
be some integers. Denote by
$\overline{m_1,m_2}$
the integer interval
$\{m_1,m_1+1,\dots,m_2\}$.
For every
$n\in{\bf N}$,
instead of
$\overline{1,n}$
we write simply
$\overline{n}$.
Also, let
$\text{\bf N}_0=\text{\bf N}\cup\{0\}$.

Let
$V$ and $\Lambda$
be finite sets, and
$V\ne\emptyset$.
By {\em unordered hypergraph\/} or simply {\em hypergraph\/} we mean an ordered pair
$H=(V,{\mathcal E})$,
where
${\mathcal E}=(e_\lambda,\lambda\in\Lambda)$
is a family of the subsets
$e_{\lambda}$
of the set $V$ (the subsets can be empty and can repeat themselves, and even the family
${\mathcal E}$
may be empty). Let us call elements of the set $V$ {\em vertices}, and members of the family
$\mathcal E$
{\em edges} of the hypergraph $H$. We write
$e_\lambda\in{\mathcal E}$ if $e_\lambda$
is a member of the family
${\mathcal E}$.
In what follows the set of vertices of a hypergraph $H$ will also be denoted by
$VH$,
and the family of its edges --- by
${\mathcal E}H$. If $|\Lambda|=m$ and $|V|=n$,
then we call hypergraph $H$ an
$(m,n)$-{\em hypergraph\/}.

Denote by
$||e||$
the {\em multiplicity\/} of an
$e\subseteq V$ in ${\mathcal E}$, i.e., $||e||=|\{\lambda\in\Lambda\,|\,e_\lambda=e\}|$.
We say that a hypergraph $H$ has {\em no\/} (or, is {\em without\/}) {\em multiple edges\/} if
$||e||\le1$
for every
$e\subseteq V$.

Note that a graph may be regarded as a special case of a hypergraph. As above, in what follows
the set of vertices of a graph $G$ is referred to as
$VG$,
and its set of edges as
$EG$.

Let
$H=(V,{\mathcal E})$, ${\mathcal E}=(e_\lambda,\lambda\in\Lambda)$,
be a hypergraph. If a linear order
$\le$
is given on the set of indices
$\Lambda$,
we say that $H$ is an {\em ordered\/} hypergraph; in that case, for the family of edges, instead of
the notation
${\mathcal E}=(e_\lambda,\lambda\in\Lambda)$
we often use the notation
${\mathcal E}=(e_\lambda,\lambda\in(\Lambda,\le))$. If $\Lambda=\{\lambda_1<\lambda_2<\dots<\lambda_m\}$,
denote by
${\mathcal E}H[i]={\mathcal E}[i]$, $i\in\overline{m}$,
the edge
$e_{\lambda_i}$
of the hypergraph $H$. If $H$ is ordered, denote by
$\overline{H}$
the corresponding unordered hypergraph.

Let
$H=(V,{\mathcal E})$, ${\mathcal E}=(e_\lambda,\lambda\in\Lambda)$,
be a hypergraph. A vertex
$v\in V$
is {\em incident to\/} an edge
$e_\lambda\in{\mathcal E}$
(or $e_\lambda$ is {\em incident to\/} $v$) if
$v\in e_\lambda$.
A vertex $v$ is called an {\em isolated vertex\/} in $H$ if there is no edge in $H$ which is incident to $v$.
A vertex $v$ is called a {\em singular vertex\/} if either
$v\in e_\lambda$
for every
$\lambda\in\Lambda$ or $v$ is an isolated vertex in $V$. A set
$V'\subseteq V$
is a set of {\em adjacent\/} vertices in $H$ if there exists an edge
$e_\lambda\in{\mathcal E}$
such that
$V'\subseteq e_\lambda$.

Let
$H=(V,{\mathcal E})$, ${\mathcal E}=(e_\lambda,\lambda\in\Lambda)$,
and
$H'=(V',{\mathcal E}')$, ${\mathcal E}'=(e'_\lambda,\lambda\in\Lambda')$,
be two hypergraphs. The hypergraph
$H'$
is a {\em subhypergraph\/} of $H$ if
$V'\subseteq V$, $\Lambda'\subseteq\Lambda$, and $e'_\lambda=e_\lambda\cap V'$
for every
$\lambda\in\Lambda'$.
If the hypergraph $H$ is ordered with liner order
$\le$
on the set
$\Lambda$,
then we suppose that the subhypergraph
$H'$
is also ordered, and that the linear order
$(\Lambda'\times\Lambda')\,\cap\le$
is given on the set
$\Lambda'$.
For every
$V_1\subseteq V$ and $\Lambda_1\subseteq\Lambda$
denote by
$H[V_1,\Lambda_1]$
the corresponding subhypergraph of $H$. So we have
$H'=H[V',\Lambda']$.

Let
$H=(V,{\mathcal E})$, ${\mathcal E}=(e_\lambda,\lambda\in\Lambda)$,
and let
$V'\subseteq V$. Take $\Lambda'=\{\lambda\in\Lambda\,|\,e_\lambda\subseteq V'\}$.
The subhypergraph
$H'=H[V',\Lambda']$
is called an {\em induced subhypergraph\/} of $H$, and we say that
$V'$
induces
$H'$;
denote by
$H[V']$
such hypergraph
$H'$.

Let
$H=(V,{\mathcal E})$, ${\mathcal E}=(e_\lambda,\lambda\in\Lambda)$,
be a hypergraph. An edge
$e_\lambda\in{\mathcal E}$
is a $k$-{\it edge\/} [$k_\le$-{\it edge\/}],
$k\in\text{\bf N}_0$, if $|e_\lambda|=k$ [$|e_\lambda|\le k$].
We also call a $0$-edge an {\em empty edge\/}, and a $|V|$-edge --- a {\em full edge\/}. If
$\emptyset\in{\mathcal E}$ [$\emptyset\not\in{\mathcal E}$],
then we say that $H$ is a hypergraph {\em with\/} [{\em without\/}] {\em empty edges\/}. We say that
$H$ is a $k_{\le}$-{\em dimensional\/} hypergraph, and write
$\text{\rm dim\,}H\le k$,
if every edge of $H$ is a $k_{\le}$-edge. We say that $H$ is a $k$-{\em dimensional\/} hypergraph, and write
$\text{\rm dim\,}H=k$, if $\text{\rm dim\,}H\le k$,
and if there exists at least one $k$-edge in $H$. The hypergraph $H$ is a $k$-{\em uniform\/} hypergraph if every
its edge is a $k$-edge. It is clear that a $k$-uniform hypergraph is always a $k$-dimensional hypergraph.

Let
$H=(V,{\mathcal E})$, ${\mathcal E}=(e_\lambda,\lambda\in\Lambda)$,
be a hypergraph. If a linear order
$\le$
is given on $V$, then we say that $H$ is a {\em hypergraph with labelling\/} ({\em determined by\/} $\le$); in that
case we also say that $H$ is a {\em hypergraph on\/}
$(V,\le)$. If $V=\{v_1<\dots<v_n\}$,
denote by
$VH[i]=V[i]$, $i\in\overline{n}$,
the vertex
$v_i$.

Let
$H_1=(V_1,{\mathcal E}_1)$, ${\mathcal E}_1=(e_\lambda,\lambda\in\Lambda_1)$,
and
$H_2=(V_2,{\mathcal E}_2)$, ${\mathcal E}_2=(e'_\lambda,\lambda\in\Lambda_2)$,
be two hypergraphs [hypergraphs with labellings]. They are {\em isomorphic\/},
$H_1\simeq H_2$ [$H_1\equiv H_2$],
if there are bijections
$\iota\!:V_1\to V_2$ and $\nu\!:\Lambda_1\to\Lambda_2$
such that
$\iota(e_\lambda)=e'_{\nu(\lambda)}$
for every
$\lambda\in\Lambda_1$ [$\iota(VH_1[i])=VH_2[i]$
for every
$i\in\overline{|V_1|}$ and $\iota(e_\lambda)=e'_{\nu(\lambda)}$
for every
$\lambda\in\Lambda_1$];
the pair
$(\iota,\nu)$
is called an {\em isomorphism\/} between
$H_1$ and $H_2$.

Let
$H_1=(V_1,{\mathcal E}_1)$, ${\mathcal E}_1=(e_\lambda,\lambda\in(\Lambda_1,\le_1))$,
and
$H_2=(V_2,{\mathcal E}_2)$, ${\mathcal E}_2=(e'_\lambda,\lambda\in(\Lambda_2,\le_2))$,
be ordered hypergraphs [ordered hypergraphs with labellings]. They are {\em isomorphic\/},
$H_1\simeq H_2$ [$H_1\equiv H_2$],
if there exists an isomorphism
$(\iota,\nu)$
between
$\overline{H}_1$ and $\overline{H}_2$
such that
$\lambda\le_1\lambda'$
iff
$\nu(\lambda)\le_2\nu(\lambda')$
for every
$\lambda,\lambda'\in\Lambda_1$.

The relations
$\simeq$ and $\equiv$
are relations of equivalence on the class of all corresponding hypergraphs. By {\em $($unlabelled$)$ hypergraph\/}
[{\em ordered $($unlabelled$)$ hypergraph\/}] we mean a class of the equivalence
$\simeq$.
By {\em labelled hypergraph\/} [{\em ordered labelled hypergraph\/}] we mean a class of the equivalence
$\equiv$.

Let $V$ be an $n$-set,
$\le$
be a linear order on $V$, and $K$ be a labelled [an ordered labelled] $(m,n)$-hypergraph. It is clear
that only one [ordered] $(m,n)$-hypergraph with labelling from the class $K$ is a hypergraph on
$(V,\le)$.
So, if we have some `hypergraph property'
$\mathfrak p$,
and want to know how many labelled [ordered labelled] $(m,n)$-hypergraphs have this property, it is
sufficient to find how many $(m,n)$-hypergraphs on
$(V,\le)$
satisfy the property
$\mathfrak p$.
This is exactly what we are going to do in the paper: when we enumerate labelled [ordered labelled] $(m,n)$-hypergraphs
satisfying the property
$\mathfrak p$,
we shall fix a linearly ordered set
$(V,\le)$
(the context will usually make it clear which
$(V,\le)$
is meant), and enumerate $(m,n)$-hypergraphs on
$(V,\le)$
satisfying the property
$\mathfrak p$.
For the sake of simplicity we shall then say `an $(m,n)$-hypergraph' instead of `an $(m,n)$-hypergraph on
$(V,\le)$'.

Let $H$ be an (unordered) hypergraph. It is a {\em multiantichain\/} if
$e_1\subseteq e_2$
implies
$e_1=e_2$
for every
$e_1,e_2\in{\mathcal E}H$
(here $e_1$ and $e_2$ taken as sets are equal). A multiantichain without multiple edges is usually called an
{\em antichain\/}. An {\em $(m,n)$-multiantichain\/} [an {\em $(m,n)$-antichain\/}] is a multiantichain [an
antichain] with $m$ edges and $n$ vertices.

Let $H$ be an (unordered) hypergraph. It is a {\em cover\/} if there is no isolated vertex in $H$. It is a
{\em proper cover\/} if it is a cover and does not contain a full edge. An {\em $(m,n)$-cover\/} [{\em a proper
$(m,n)$-cover\/}] is a cover [a proper cover] with $m$ edges and $n$ vertices.

Let us agree that in the definition that follows the symbol ${<}\beta{>}$ means only one of the following words:
`hypergraph' or `multiantichain' or `antichain' or `cover' or one of these words with the prefix `$(m,n)$\,-'. By
analogy with the notions of $T_0$-, $T_1$- and $T_2$-spaces from general topology let us introduce similar notions
for hypergraphs. A ${<}\beta{>}$ $H$ is a:
\begin{enumerate}
\item[a)]
$T_0$-${<}\beta{>}$
iff for every two different vertices
$u,v\in VH$
there exists
$e\in{\mathcal E}H$
such that
$(u\in e\land v\not\in e)\lor (u\not\in e\land v\in e)$,
\item[b)]
$T_1$-${<}\beta{>}$
iff for every pair
$(u,v)\in (VH)^2$, $u\ne v$,
there exists
$e\in{\mathcal E}H$
such that
$(u\in e\land v\not\in e)$,
\item[c)]
$T_2$-${<}\beta{>}$
iff for every pair
$(u,v)\in (VH)^2$, $u\ne v$,
there exist two edges
$e_1,e_2\in{\mathcal E}H$
such that
$(u\in e_1\land v\in e_2\land e_1\cap e_2=\emptyset)$.
\end{enumerate}

In the case of a proper cover we speak about proper $T_i$-cover and proper
$T_i$-$(m,n)$-cover, $i\in\overline{0,2}$.
Note that every
$T_i$-hypergraph,
$i\in\overline{2}$,
with at least two vertices, is a cover. It is also clear that every $T_0$-hypergraph can have at most
one isolated vertex.

Let $H$ be a hypergraph, and let
$u,v\in VH$
be its two vertices (the vertices can be equal or different). A {\em path\/} in $H$ connecting the
vertices $u$ and $v$ is either an edge
$e\in{\mathcal E}H$
such that
$u,v\in e$,
or a finite sequence
$e_1,e_2,\dots,e_k$
of $k$ different edges of $H$,
$k\ge 2$,
such that
$u\in e_1$, $v\in e_k$ and $e_i\cap e_{i+1}\ne\emptyset$
for every
$i\in\overline{1,k-1}$.
The hypergraph $H$ is {\em connected\/} if for every its two vertices there exists a path connecting them.
So every connected hypergraph is a cover. Let
$V_1\subseteq V$,
and suppose that for every
$\lambda\in\Lambda$,
either
$e_{\lambda}\cap V_1=e_{\lambda}$ or $e_{\lambda}\cap V_1=\emptyset$.
If the hypergraph
$H[V_1]$
is connected, then
$V_1$
is called a {\em component\/} ({\em component of connectedness\/}) of $H$; if
$|V_1|=k$,
we also say that
$V_1$
is a {\em $k$-component\/}.

A hypergraph
$H=(V,{\mathcal E})$, ${\mathcal E}=(e_\lambda,\lambda\in\Lambda)$,
has $k$-{\em intersecting\/} property,
$k\ge2$, if $e_{\lambda_1}\cap\dots\cap e_{\lambda_k}\ne\emptyset$
for every $k$ edges
$e_{\lambda_1},\dots,e_{\lambda_k}\in{\mathcal E}$ (some of $\lambda_i$ may be equal).
If the hypergraph $H$ has $k$-intersecting property for every
$k\ge2$, i.e., $\cap_{\lambda\in\Lambda}e_\lambda\ne\emptyset$,
we say that the hypergraph has {\em intersecting property\/}, or, for short, that it has {\em $\cap$-property\/};
if a hypergraph has no intersecting property, sometimes we say that it has {\em $\neg\cap$-property\/}. If $H$
has $\cap$-property, then it has at least one singular vertex.

Each of the introduced hypergraphs may be ordered or unordered, with a labelling or without it. Each of the
properties determining such hypergraphs is invariant with respect to
$\equiv$ and $\cong$,
and consequently we can speak about corresponding labelled and unlabelled hypergraphs. In the following we are
dealing only with labelled hypergraphs, and, from now on, we leave out the word `labelled' and keep in mind that
all hypergraphs are labelled. Note that all problems investigated in the paper can also be formulated for
unlabelled case but we cannot solve most of them yet.

Let us agree that we denote by ${<}\alpha{>}$ the word combination `unordered [ordered]'.

Let
$\mathfrak{p}$
be a property that an ${<}\alpha{>}$ hypergraph could have, or, as we often say, an ${<}\alpha{>}$ {\em hypergraph
property\/}. If an ${<}\alpha{>}$ hypergraph $H$ has the property
$\mathfrak{p}$,
we say also that $H$ is an ${<}\alpha{>}$ {\em $\mathfrak{p}$-hypergraph\/}. Now let
$\mathfrak{p}_1,\dots,\mathfrak{p}_k$
be some ${<}\alpha{>}$ hypergraph properties. If an ${<}\alpha{>}$ hypergraph has all these properties, we say
that it has the property
$\mathfrak{p}_1\land\dots\land\mathfrak{p}_k$,
and, consequently, it is an
${<}\alpha{>}$ $\mathfrak{p}_1\land\dots\land\mathfrak{p}_k$-hypergraph,
or, as we sometimes say, an
${<}\alpha{>}$ $\mathfrak{p}_1$-\dots-$\mathfrak{p}_k$-hypergraph.

Denote by
$\mathfrak{H}_{\forall}$ [$\vec{\mathfrak{H}}_{\forall}$]
the class of all ${<}\alpha{>}$ hypergraphs. Let
$\mathfrak{p}$
be an ${<}\alpha{>}$ hypergraph property. Denote by
$\mathfrak{H}_{\mathfrak{p}}$ [$\vec{\mathfrak{H}}_{\mathfrak{p}}$]
the class of all ${<}\alpha{>}$ $\mathfrak{p}$-hypergraphs, and by
$\mathfrak{H}^\ast_{\mathfrak{p}}$ [$\vec{\mathfrak{H}}\vphantom{\mathfrak{H}}^\ast_{\mathfrak{p}}$]
the class of all ${<}\alpha{>}$ $T_0$-$\mathfrak{p}$-hypergraphs. If
$\mathfrak{p}=\mathfrak{p}_1\land\dots\land\mathfrak{p}_k$,
then
$$
\mathfrak{H}_{\mathfrak{p}}=\mathfrak{H}_{\mathfrak{p}_1\land\dots\land\mathfrak{p}_k}=
\mathfrak{H}_{\mathfrak{p}_1}\cap\dots\cap\mathfrak{H}_{\mathfrak{p}_k}
\quad[\vec{\mathfrak{H}}_{\mathfrak{p}}=\vec{\mathfrak{H}}_{\mathfrak{p}_1\land\dots\land\mathfrak{p}_k}=
\vec{\mathfrak{H}}_{\mathfrak{p}_1}\cap\dots\cap\vec{\mathfrak{H}}_{\mathfrak{p}_k}];
$$
so,
$\mathfrak{H}^\ast_{\mathfrak{p}}=\mathfrak{H}_{T_0\land\mathfrak{p}}$
[$\vec{\mathfrak{H}}\vphantom{\mathfrak{H}}^\ast_{\mathfrak{p}}=\vec{\mathfrak{H}}_{T_0\land\mathfrak{p}}$].
If
$\mathfrak{H}$
is a class of ${<}\alpha{>}$ hypergraphs, then denote by
$\mathfrak{H}(m,n)$
the class of all $(m,n)$-hypergraphs from
$\mathfrak{H}$.
Take
$\alpha_{\mathfrak{p}}(m,n)=|\mathfrak{H}_{\mathfrak{p}}(m,n)|$
[\,$\vec{\alpha}_{\mathfrak{p}}(m,n)=|\vec{\mathfrak{H}}_{\mathfrak{p}}(m,n)|$\,]
and
$\alpha^\ast_{\mathfrak{p}}(m,n)=|\mathfrak{H}_{\mathfrak{p}}^\ast(m,n)|$
[\,$\vec{\alpha}\vphantom{\alpha}^\ast_{\mathfrak{p}}(m,n)=
    |\vec{\mathfrak{H}}\vphantom{\mathfrak{H}}^\ast_{\mathfrak{p}}(m,n)|$\,].

The {\it incidence matrix\/} of a
$H\in\vec{\mathfrak{H}}_{\forall}(m,n)$
is the binary matrix
$M_H=[m_{ij}]_{m\times n}$,
where for every
$i\in\overline{m}$ and $j\in\overline{n}$,
$m_{ij}=1$ if $V[j]\in{\mathcal E}[i]$, and $m_{ij}=0$ if $V[j]\not\in{\mathcal E}[i]$.
Introducing in such a way the incidence matrix for an ordered $(m,n)$-hypergraph we define a bijective map
between the class
$\vec{\mathfrak{H}}_{\forall}(m,n)$
and the class of all binary matrices with $m$ rows and $n$ columns. So, if we enumerate a class of ordered
$\mathfrak{p}$-hypergraphs [$T_0$-$\mathfrak{p}$-hypergraphs],
we in fact enumerate the corresponding class of binary matrices [binary matrices with different columns].

The {\em dual hypergraph\/}
$H^T$
of an ordered labelled hypergraph
$H\in\vec{\mathfrak{H}}_{\forall}(m,n)$
is a hypergraph from
$\vec{\mathfrak{H}}_{\forall}(n,m)$
whose incidence matrix is
$M_H^T$,
where
$M_H^T$
is the transpose of
$M_H$.
Let an ordered hypergraph property
$\mathfrak{p}$
be given.
We say that a hypergraph
$H\in\vec{\mathfrak{H}}_{\forall}(m,n)$
has the {\em dual property\/} of
$\mathfrak{p}$,
that is, $H$ is a {\em dual $\mathfrak{p}$-hypergraph\/}, or, as we also say, $H$ is {\it dually\/}
$\mathfrak{p}$,
if $H^T$ has the property $\mathfrak{p}$.

\section{On $T_0$-hypergraphs}

As we have agreed, from now on every hypergraph will be labelled, and instead of an ordered [an unordered] labelled
hypergraph we simply say an ordered [an unordered] hypergraph.

Let
$H=(V,{\mathcal E})$, ${\mathcal E}=(e_\lambda,\lambda\in\Lambda)$,
be an arbitrary hypergraph (ordered or unordered). Let us define a relation
$\sim_H$
on $V$ such that for every
$u,v\in V$,
$$
u\sim_H v \Leftrightarrow(\forall e\in{\mathcal E})(u\in e\land v\in e)\lor(u\not\in e\land v\not\in e).
$$
It is easy to see that the following proposition holds.

\begin{proposition}\label{pp31}
  For every ordered $[$unordered$\,]$ hypergraph $H$, the relation
  $\sim_H$
  is a relation of equivalence.
\end{proposition}

Let
$H=(V,{\mathcal E})$, ${\mathcal E}=(e_\lambda,\lambda\in\Lambda)$,
be an arbitrary hypergraph. For every
$v\in V$,
denote by
$[v]$
the class of equivalence
$\sim_H$
containing $v$. Also, denote by
$\pi_H$
the partition of the set $V$ corresponding to the equivalence
$\sim_H$, i.e.\ $\pi_H=V/{\sim_H}$.
It is easy to see that for every
$e\in{\mathcal E}$ and $v\in V$, if $e\cap[v]\ne\emptyset$,
then
$[v]\subseteq e$.
Let
$$
e'_\lambda=[e_\lambda]=\{[v]\in \pi_H\,|\,[v]\subseteq e_\lambda\}\qquad\text{for every}\;\lambda\in\Lambda,
$$
and let
$[{\mathcal E}]=(e'_\lambda\,|\,\lambda\in\Lambda\}$.
Denote by
$[H]$
the hypergraph
$(V/\sim_H,[{\mathcal E}])$.

\begin{proposition}\label{pp32}
  For every ordered $[$unordered$\,]$ hypergraph
  $H=(V,{\mathcal E})$,
  the hypergraph
  $[H]$
  is a
  $T_0$-hypergraph.
\end{proposition}

{\it Proof.\/}\;
Let
$[u]$ and $[v]$, $[u]\ne[v]$,
be two different elements of
$\pi_H$.
As
$[u]\ne[v]$,
then
$\neg(u\sim_H v)$,
and therefore there exists
$e\in{\mathcal E}$
such that
$(u\in e\land v\not\in e)\lor(u\not\in e\land v\in e)$.
Consequently, we get that
$([u]\in [e]\land [v]\not\in [e])\lor([u]\not\in [e]\land [v]\in [e])$,
and, therefore, the hypergraph
$[H]$
is a
$T_0$-hypergraph.\,$\square$

By a {\em partition type\/} we mean any $n$-tuple
$\tau=(\alpha_1,\dots,\alpha_n)$
such that
$\sigma(\tau)=\alpha_1+2\alpha_2+\dots+n\alpha_n=n$, and $\alpha_i\in\text{\bf N}_0$
for every
$i\in\overline{n}$;
denote by
$|\tau|$
the number
$\alpha_1+\alpha_2+\dots+\alpha_n$.

Let
$\tau=(\alpha_1,\dots,\alpha_n)$
be a partition type. We say that a partition
$\pi$
of an $n$-set has the type
$\tau$
if it has
$\alpha_i$
partition classes of the cardinality $i$ for every
$i\in\overline{n}$.
The number of all partition classes of
$\pi$
is denoted by
$|\pi|$,
and the partition type
$\tau$
of the partition
$\pi$
is denoted by
$\text{\rm typ}(\pi)$.
It is clear that
$|\pi|=|\text{\rm typ}(\pi)|$.

Let
$\tau=(\alpha_1,\dots,\alpha_n)$
be a partition type. Denote by
$b(\tau)$
the number of all partitions
$\pi$
of a given $n$-set such that
$\text{\rm typ}(\pi)=\tau$.
It is well known that
$$
b(\tau)=\dfrac{n!}{\alpha_1!\alpha_2!\dots\alpha_n!(2!)^{\alpha_2}\dots(n!)^{\alpha_n}}.
$$

Fix a countable set
$W_\infty=\{w_i\,|\,w_i\in\text{\bf N}_0\}$,
and for every
$n\in\text{\bf N}$,
let
$W_n=\{w_1<w_2<\ldots<w_n\}$.
Denote by
${\mathfrak G}_n$
the set of all labelled graphs on
$W_n$.

Denote by
$\Pi(n)$
the set of all partitions of the $n$-set
$W_n$,
and denote by
$\Pi(n,i)$
the set of all partitions of
$W_n$
into $i$ parts. Let
$\tau=(\alpha_1,\dots,\alpha_n)$
be a partition type. Denote by
$\Pi(\tau)=\Pi(\alpha_1,\dots,\alpha_n)$
the set of all partitions of the type
$\tau$
from
$\Pi(n)$.

Let
$\mathfrak{p}$
be an unordered hypergraph property. An $(m,n)$-$\mathfrak{p}$-hypergraph $H$ satisfies the property
$p_{ij}$ if $V[i]\sim_H V[j]$; $i,j\in\overline{n}$, $i\ne j$.
Let
${\mathcal P}=\{p_{i_1j_1},\dots,p_{i_kj_k}\}$
be a subset of
${\mathcal P}_0=\{p_{ij}\,|\,i,j\in\overline{n}, i\ne j\}$.
Denote by
$\alpha_{\mathfrak{p}}(m,n;{\mathcal P})$
the number of all $\mathfrak{p}$-hypergraphs satisfying at least all the properties from
${\mathcal P}$.
Note that
$\alpha_{\mathfrak{p}}(m,n;\emptyset)=\alpha_{\mathfrak{p}}(m,n)$.
Consider the graph
$G({\mathcal P})=(W_n,E({\mathcal P}))$,
where
$E({\mathcal P})=\{\{w_{i_s},w_{j_s}\}\,|\,s\in\overline{k}\}$.
Find the components of
$G({\mathcal P})$,
and denote by
$\pi({\mathcal P})$
the corresponding partition of the set
$W_n$. Let ${\mathcal P}$ and ${\mathcal P}'$
be subsets of
${\mathcal P}_0$
such that
$\pi({\mathcal P})=\pi({\mathcal P}')$.
It is clear that every $\mathfrak{p}$-hypergraph $H$ which satisfies all the properties from
${\mathcal P}$
satisfies all the properties from
${\mathcal P}'$,
and vice versa. Consequently, we get the following lemma.

\begin{lemma}\label{lmm31}
  For every
  ${\mathcal P},{\mathcal P}'\subseteq{\mathcal P}_0$,
  if
  $\pi({\mathcal P})=\pi({\mathcal P}')$,
  then
  $\alpha_{\mathfrak{p}}(m,n;{\mathcal P})=\alpha_{\mathfrak{p}}(m,n;{\mathcal P}')$.
\end{lemma}

Denote by
$\alpha_{\mathfrak{p}}(m,n;\pi)$
the value
$\alpha_{\mathfrak{p}}(m,n;{\mathcal P})$,
where
${\mathcal P}\subseteq{\mathcal P}_0$
such that
$\pi({\mathcal P})=\pi$.
From the above lemma it follows that this notation is well founded. Now by using the above lemma
we can prove the following theorem.

\begin{theorem}\label{trm1}
  $$
  \alpha^\ast_{\mathfrak{p}}(m,n)=\displaystyle
  \sum_{(\alpha_1,\dots,\alpha_n)}\,
  \sum_{\pi\in\Pi(\alpha_1,\dots,\alpha_n)}\alpha_{\mathfrak{p}}(m,n;\pi)
  \left[\prod_{i=1}^n[(-1)^{i-1}(i-1)!]^{\alpha_i}\right].
  $$
\end{theorem}

{\it Proof.}\;
Using the lemma and the formula of inclusion and exclusion we get
$$
\begin{array}{lcl}
  \alpha^\ast_{\mathfrak{p}}(m,n)
 & = &
  \displaystyle\sum_{{\mathcal P}\subseteq{\mathcal P}_0}(-1)^{|{\mathcal P}|}\alpha_{\mathfrak{p}}(m,n;{\mathcal P})=
  \sum_{{\mathcal P}\subseteq{\mathcal P}_0}(-1)^{|EG({\mathcal P})|}\alpha_{\mathfrak{p}}(m,n;{\mathcal P})=\\[6pt]
 & = &
  \displaystyle\sum_{\pi\in\Pi(n)}\sum_{{\mathcal P}\subseteq{\mathcal P}_0,\,\pi({\mathcal P})=\pi}(-1)^{|EG({\mathcal P})|}
  \alpha_{\mathfrak{p}}(m,n;{\mathcal P})=\\[6pt]
 & = &
  \displaystyle\sum_{\pi\in\Pi(n)}\mu(\pi)\,\alpha_{\mathfrak{p}}(m,n;\pi)=\sum_{(\alpha_1,\dots,\alpha_n)}\,
  \sum_{\pi\in\Pi(\alpha_1,\dots,\alpha_n)}\mu(\pi)\,\alpha_{\mathfrak{p}}(m,n;\pi),\\
\end{array}
$$
where
$$
\mu(\pi)=\sum_{{\mathcal P}\subseteq{\mathcal P}_0,\,\pi({\mathcal P})=\pi}(-1)^{|EG({\mathcal P})|}.
$$

Let
$(\alpha_1,\dots,\alpha_n)$
be a partition type of an $n$-set, and let
$\pi=\{V_1,\dots,V_k\}\in\Pi(\alpha_1,\dots,\alpha_n)$.
Take a
${\mathcal P}\subseteq{\mathcal P}_0$
such that
$\pi({\mathcal P})=\pi$.
Denote by
$\mathfrak{G}(\pi)$
the class of all graphs from
$\mathfrak{G}_n$
having $k$ components of connectedness defined by the sets
$V_1,\dots,V_k$.
Also for every set $V$, denote by
$\bar{\mathfrak{G}}(V)$
the set of all connected graphs on $V$ (having $V$ as the set of their vertices). Then it is clear that
$$
\mu(\pi)=\sum_{G\in\mathfrak{G}(\pi)}(-1)^{|EG|}=
\sum_{G_1\in\bar{\mathfrak{G}}(V_1)}\ldots\sum_{G_k\in\bar{\mathfrak{G}}(V_k)}(-1)^{|EG_1|}\dots(-1)^{|EG_k|}.
$$
As for every $n$-set $V$
$$
\sum_{G\in\bar{\mathfrak{G}}(V)}(-1)^{|EG|}=(-1)^{n-1}(n-1)!,
$$
then
$$
\begin{array}{lcl}
   \mu(\pi)
   & = &
   \displaystyle\left[\sum_{G_1\in\bar{\mathfrak{G}}(V_1)}(-1)^{|EG_1|}\right]\ldots
   \left[\sum_{G_k\in\bar{\mathfrak{G}}(V_k)}\dots(-1)^{|EG_k|}\right]=\\[6pt]
   & = &
   \displaystyle\prod_{i=1}^n[(-1)^{i-1}(i-1)!]^{\alpha_i}.\,\square
\end{array}
$$

\smallskip
We say that
$\mathfrak{p}$
is a {\em uniform property\/} if the equality
$\text{\rm typ}(\pi)=\text{\rm typ}(\pi')=\tau$
always implies the equality
$\alpha_{\mathfrak{p}}(m,n;\pi)=\alpha_{\mathfrak{p}}(m,n;\pi')=\alpha_{\mathfrak{p}}(m,n;\tau)$.

Let
$\tau$
be a partition type of an $n$-set. Let us put
$$
c(\tau)=\dfrac{n!}{\alpha_1!\alpha_2!\dots\alpha_n!2^{\alpha_2}\dots n^{\alpha_n}}.
$$
We have the following theorem.

\begin{theorem}\label{trm2}
  If
  $\mathfrak{p}$
  is a uniform property, then
  $$
  \alpha^\ast_{\mathfrak{p}}(m,n)=\sum_{\tau=(\alpha_1,\dots,\alpha_n),\,\sigma(\tau)=n}(-1)^{n-\alpha_1-\dots-\alpha_n}
  c(\tau)\,\alpha_{\mathfrak{p}}(m,n;\tau),
  $$
\end{theorem}

{\it Proof.}\;From Theorem \ref{trm1} it follows that
$$
\begin{array}{lcl}
 \alpha^\ast_{\mathfrak{p}}(m,n) & = & \displaystyle
 \sum_{(\alpha_1,\dots,\alpha_n)}\,
 \sum_{\pi\in\Pi(\alpha_1,\dots,\alpha_n)}\alpha_{\mathfrak{p}}(m,n;\pi)
 \left[\prod_{i=1}^n[(-1)^{i-1}(i-1)!]^{\alpha_i}\right]=\\
                                            & = & \displaystyle
 \sum_{\tau=(\alpha_1,\dots,\alpha_n),\,\sigma(\tau)=n}\,
 \alpha_{\mathfrak{p}}(m,n;\tau)\,b(\tau)
 \left[\prod_{i=1}^n[(-1)^{i-1}(i-1)!]^{\alpha_i}\right],
\end{array}
$$
and we immediately get the formula.\,$\square$

\smallskip
We say that a uniform property
$\mathfrak{p}$
is a {\em regular $T_0$-property\/} if the equality
$|\tau|=|\tau'|=k$
always implies the equality
$\alpha_{\mathfrak{p}}(m,n;\tau)=\alpha_{\mathfrak{p}}(m,n;\tau')=\alpha_{\mathfrak{p}}(m,n;k)$.
Now we have the following theorem.

\begin{theorem}\label{trm3}
  If
  $\mathfrak{p}$
  is a regular $T_0$-property, then
  \begin{equation*}
    \alpha_{\mathfrak{p}}^\ast(m,n)=\sum^{n}_{i=1}s_{n,i}\,\alpha_{\mathfrak{p}}(m,n,i),
  \end{equation*}
  where
  $s_{n,i}$
  are the Stirling numbers of the first kind.
\end{theorem}

{\it Proof.}\; From Theorem \ref{trm2} it follows that
$$
\begin{array}{lcl}
 \alpha^\ast_{\mathfrak{p}}(m,n)
& = & \displaystyle
 \sum_{\tau=(\alpha_1,\dots,\alpha_n),\,\sigma(\tau)=n}(-1)^{n-\alpha_1-\dots-\alpha_n}
 c(\tau)\,\alpha_{\mathfrak{p}}(m,n;\tau)=\\
& = & \displaystyle
 \sum_{k=1}^{n}\alpha_{\mathfrak{p}}(m,n;k)\sum_{\tau=(\alpha_1,\dots,\alpha_n),\,|\tau|=k,\,\sigma(\tau)=n}\,
 (-1)^{n-k}c(\tau).
\end{array}
$$
Now from identity (see \cite{Rio})
$$
\sum_{\tau=(\alpha_1,\dots,\alpha_n),\,|\tau|=k}\,(-1)^{n-k}c(\tau)=s_{n,k}
$$
follows the theorem.\,$\square$

\smallskip
Let $V$ be an $n$-set, and let
$V_i$, $i\in\overline{k}$,
be partition classes of a $k$-partition
$\pi$
of $V$. An ordered [unordered] hypergraph
$H=(V,{\mathcal E})$, ${\mathcal E}=(e_\lambda,\lambda\in\Lambda)$,
is {\em $\pi$-granular\/} if either
$e_\lambda\cap V_i=\emptyset$ or $e_\lambda\cap V_i=V_i$
for every
$i\in\overline{k}$ and $\lambda\in\Lambda$.
For example, every hypergraph $H$ is $\pi_H$-granular.

Let $V$ be a linearly ordered $n$-set, let
$\pi=(V_1,\dots,V_k)$
be an ordered $k$-partition,
$k\le n$, of $V$,
and let
$H=(V,{\mathcal E})$, ${\mathcal E}=(e_\lambda,\lambda\in\Lambda)$,
be a $\pi$-granular (labelled) ordered [unordered] $(m,n)$-hypergraph on $V$. An ordered [unordered]
$(m,k)$-hypergraph
$H'=(V',{\mathcal E}')$, ${\mathcal E}'=(e'_\lambda,\lambda\in\Lambda)$,
is a {\em $\pi$-condensation\/} of $H$ (or, $H$ is a {\em $\pi$-expansion\/} of
$H'$) if $V'[i]\in e'_\lambda$ iff $V_i\subseteq e_\lambda$.
We say that
$H'$
is a condensation of $H$ (or, $H$ is an expansion of
$H'$)
if there exists a partition
$\pi$ of $VH$
such that
$H'$
is a $\pi$-condensation of $H$ (or, $H$ is a $\pi$-expansion of
$H'$).

We say that a regular $T_0$-property
$\mathfrak{p}$
{\em allows $T_0$-filtration\/} if
$\alpha_{\mathfrak{p}}(m,n;k)=\alpha_{\mathfrak{p}}(m,k)$
for every
$k\in\overline{n}$.
It is not difficult to see that a regular $T_0$-property
$\mathfrak{p}$
allows $T_0$-filtration if it is invariant to the operations of condensation and expansion, i.e., if a
hypergraph $H$ has the property
$\mathfrak{p}$,
then every expansion [condensation] of $H$ has the property
$\mathfrak{p}$.
Now we are able to generalize Lemma 1 from \cite{JK1}. From Theorem \ref{trm3} we get

\begin{theorem}\label{trm4}
  Let
  $\mathfrak{p}$
  be a property that allows $T_0$-filtration. Then
  \begin{equation}
    \alpha_{\mathfrak{p}}^\ast(m,n)=\sum^{n}_{i=1}s_{n,i}\,\alpha_{\mathfrak{p}}(m,i).\label{f1}
  \end{equation}
\end{theorem}

\begin{theorem}\label{trm5}
  Let
  $\mathfrak{p}$
  be a property that allows $T_0$-filtration. Then
  \begin{equation}
    \alpha_{\mathfrak{p}}(m,n)=\sum^{n}_{i=1}S_{n,i}\,\alpha^\ast_{\mathfrak{p}}(m,i),\label{f2}
  \end{equation}
  where
  $S_{n,i}$
  are the Stirling numbers of the second kind.
\end{theorem}

{\it Proof.\/}\;
The formula can be obtained by applying the Stirling inversion (see \cite{Agn}) on (\ref{f1}).\,$\square$

\smallskip
If
$\mathfrak{p}$
is a property that allows $T_0$-filtration, and we have that the equality from Theorem~\ref{trm4} holds,
we say that {\em $F_0$-transformation\/} can be applied on
$\mathfrak{p}$,
and we write
${\mathfrak H}^\ast_{\mathfrak{p}}=F_0[{\mathfrak H}_{\mathfrak{p}}]$.
{\em $F_0^{-1}$-transformation\/} is defined by formula (\ref{f2}), and we write
${\mathfrak H}_{\mathfrak{p}}=F_0^{-1}[{\mathfrak H}^\ast_{\mathfrak{p}}]$.
$F_0$-transformation [$F_0^{-1}$-transformation] should be understood as a rule, by the application of which
every (or almost every) number
$\alpha^\ast_{\mathfrak{p}}(m,n)$ [$\alpha_{\mathfrak{p}}(m,n)$]
is obtained from a set of numbers
$\alpha_{\mathfrak{p}}(i,j)$ [$\alpha^\ast_{\mathfrak{p}}(i,j)$],
and the rule is defined by the formula (\ref{f1}) [(\ref{f2})].

Let
$\mathfrak{p}$
be a hypergraph property. Denote by
$\hat{\mathfrak{H}}_{\mathfrak{p}}(n)$ [$\hat{\mathfrak{H}}^\ast_{\mathfrak{p}}(n)$]
the class of all hypergraphs with $n$ vertices and without multiple edges
from
$\mathfrak{H}_{\mathfrak{p}}$ [$\mathfrak{H}^\ast_{\mathfrak{p}}$].
Let
$\hat{\alpha}_{\mathfrak{p}}(n)=|\hat{\mathfrak{H}}_{\mathfrak{p}}(n)|$
and
$\hat{\alpha}^\ast_{\mathfrak{p}}(n)=|\hat{\mathfrak{H}}_{\mathfrak{p}}^\ast(n)|$.
As
$\alpha_{\mathfrak{p}}(m,n)=0$ and $\alpha^\ast_{\mathfrak{p}}(m,n)=0$
for every
$m>2^n$, then
$$
\hat\alpha_{\mathfrak{p}}(n)=\sum_{i=0}^{2^n}\alpha_{\mathfrak{p}}(i,n)\quad\text{\rm and}\quad
\hat\alpha^\ast_{\mathfrak{p}}(n)=\sum_{i=0}^{2^n}\alpha^\ast_{\mathfrak{p}}(i,n).
$$
It is obvious that
$\alpha^\ast_{\mathfrak{p}}(0,n)=0$ if $n\ge2$.

\begin{theorem}\label{trm6}
  Let
  $\mathfrak{p}$
  be a property that allows
  $T_0$-filtration.
  Then
  $$
  \hat\alpha^\ast_{\mathfrak{p}}(n)=\sum^{n}_{i=0}s_{n,i}\,\hat\alpha_{\mathfrak{p}}(i).
  $$
\end{theorem}

{\em Proof.}\;
From Theorem \ref{trm4} we immediately get
$$
\hat\alpha^\ast_{\mathfrak{p}}(n)=\sum_{i=0}^{2^n}\alpha^\ast_{\mathfrak{p}}(i,n)=
\sum_{i=0}^{2^n}\sum_{j=1}^ns_{n,j}\alpha_{\mathfrak{p}}(i,j)=
\sum_{j=0}^ns_{n,j}\sum_{i=1}^{2^n}\alpha_{\mathfrak{p}}(i,j)=
\sum_{j=0}^ns_{n,j}\hat\alpha_{\mathfrak{p}}(i).\,\square
$$

\smallskip
Most of the classes of hypergraphs considered in the paper are obtained as an intersection of some
basic classes of hypergraphs. Because of that we shall often make use of the following proposition
which is not difficult to prove.

\begin{proposition}\label{pp33}
  Let
  $\mathfrak{p}_1$ and $\mathfrak{p}_2$
  be two properties which allow
  $T_0$-filtration.
  Then property
  $\mathfrak{p}_1\land\mathfrak{p}_2$
  also allows
  $T_0$-filtration.
\end{proposition}

Let us give by the following theorem the first set of properties allowing $T_0$-filtration. We can see, and
it will be clear later, that most of the properties are basic.

\begin{theorem}\label{trm7}
  Each of the following properties allows $T_0$-filtration: `to be an unordered$/$\!or\-dered hypergraph',
  `to be without$/$\!with empty edges', `to be without$/$\!with isolated vertices', `to be without$/$\!with
  full edges', `to have$/$\!not to have intersecting property', `to be without$/$\!with multiple edges', and
  `to be without$/$\!with singular vertices'.
\end{theorem}

Note that the property `to be without isolated vertices' is equivalent to the property `to be a cover'; the
dual property of the property `to be without [with] isolated vertices' is `to be without [with] empty edges';
the dual property of the property `to have [not to have] intersecting property' is `to be without [with] full
edges'; `to be with singular vertices' means that the hypergraph has intersection property or has an isolated
vertex. Also, `to be without singular vertices' means `to be a cover without intersecting property'.

Let
$H=(V,{\mathcal E})$, ${\mathcal E}=(e_\lambda,\lambda\in\Lambda)$,
and
$H'=(V,{\mathcal E}')$, ${\mathcal E}'=(e'_\delta,\delta\in\Delta)$,
be two unordered hypergraphs. Suppose that there exist a partition
$\pi=\{\Lambda_\delta\,|\,\delta\in\Delta\}$
of the set
$\Lambda$
such that
$e_\lambda=e'_\delta$ if $\lambda\in\Lambda_\delta$.
Then we say that
$H'$
is an {\em edge $\pi$-contraction\/} of $H$, and $H$ is an {\em edge $\pi$-expansion\/} of
$H'$.
An unordered hypergraph property
$\mathfrak{p}$
is {\em edge partition invariant\/} if
$\mathfrak{p}$
is invariant with respect to the operations of edge contractions and edge expansions.

\begin{theorem}\label{trm8}
  Let
  ${\mathfrak p}$
  be an edge partition invariant unordered hypergraph property. Denote by
  $\tilde{\mathfrak H}_{\mathfrak p}(m,n)$
  the class of all unordered hypergraphs without multiple edges from the class
  ${\mathfrak H}_{\mathfrak p}(m,n)$,
  and let
  $\tilde{\alpha}_{\mathfrak p}(m,n)=|\tilde{\mathfrak H}_{\mathfrak p}(m,n)|$.
  Then
  $$
  \alpha(m,n)=\sum^{m}_{i=1}C^{i-1}_{m-1}\,\tilde{\alpha}(i,n).
  $$
\end{theorem}

Let
$H=(V,{\mathcal E})$, ${\mathcal E}=(e_\lambda,\lambda\in(\Lambda,\le_1))$,
and
$H'=(V,{\mathcal E}')$, ${\mathcal E}'=(e'_\delta,\delta\in(\Delta,\le_2))$,
be two ordered hypergraphs. Suppose that there exists an ordered partition
$\pi=(\Lambda_\delta,\delta\in(\Delta,\le_2))$
of the set
$\Lambda$
such that
$e_\lambda=e'_\delta$ if $\lambda\in\Lambda_\delta$.
Then we say that
$H'$
is an {\em edge $\pi$-contraction\/} of $H$, and $H$ is an {\em edge $\pi$-expansion\/} of
$H'$.
An ordered hypergraph property
$\mathfrak{p}$
is {\em edge partition invariant\/} if
$\mathfrak{p}$
is invariant with respect to the operations of edge contractions and edge expansions.

\begin{theorem}\label{trm9}
  Let
  ${\mathfrak p}$
  be an edge partition invariant ordered hypergraph property. Denote by
  $\tilde{\mathfrak H}_{\mathfrak p}(m,n)$
  the class of all ordered hypergraphs without multiple edges from the class
  ${\mathfrak H}_{\mathfrak p}(m,n)$,
  and let
  $\tilde{\alpha}_{\mathfrak p}(m,n)=|\tilde{\mathfrak H}_{\mathfrak p}(m,n)|$.
  Then
  $$
  \alpha(m,n)=\sum^{m}_{i=1}S_{m,i}\,\tilde{\alpha}(i,n).
  $$
\end{theorem}

If
$\mathfrak{p}$
is an edge partition invariant ordered [unordered] hypergraph property, and the equality from the Theorem
\ref{trm9} [\ref{trm8}] holds, then we say that the {\em $G_1$-transformation\/} [{\em $G'_3$-transformation\/}]
can be applied on
$\mathfrak{p}$,
and we write
${\mathfrak H}_{\mathfrak{p}}=\text{\rm G}_1(\tilde{\mathfrak H}_{\mathfrak{p}})$
[${\mathfrak H}_{\mathfrak{p}}=\text{\rm G}'_3(\tilde{\mathfrak H}_{\mathfrak{p}})$].
The formulas from Theorems \ref{trm8} and \ref{trm9} are easily inverted, and they define the {\em $G_1^{-1}$-transformation\/}
[{\em $(G_3')^{-1}$-transformation\/}]. At the same time it is necessary to note that Theorems \ref{trm5} and \ref{trm9}
are identical in principle. It is just that the former concerns vertices, and the latter --- edges of a
hypergraph.

Generally speaking, our goal is to find, for a given ordered [unordered] hypergraph property
$\mathfrak{p}$,
the number of all ordered [unordered] $T_0$-$(m,n)$-$\mathfrak{p}$-hyper\-graphs. The denotations of the numbers
that we are trying to calculate in the context of the problem will always be given by the table of the
following type:
{\small\begin{center}
\begin{tabular}{|l|l|l|l|l|} \hline
            & $\mathfrak{q}_0$  & $\mathfrak{q}_1$  & .\;\,.\;\,.\;\,.\;\,.\;\,.\;\, & $\mathfrak{q}_k$\\
\hline
\parbox{1.6in}{ordered $\mathfrak{q}$-hypergraph\newline without multiple edges}
   & \parbox{0.7in}{$\alpha_{01}(m,n)$\newline $\alpha_{01}^\ast(m,n)$}
   & \parbox{0.7in}{$\alpha_{11}(m,n)$\newline $\alpha_{11}^\ast(m,n)$}
   & .\;\,.\;\,.\;\,.\;\,.\;\,.\;\,
   & \parbox{0.7in}{$\alpha_{k1}(m,n)$\newline $\alpha_{k1}^\ast(m,n)$} \\
\hline
\parbox{1.6in}{ordered $\mathfrak{q}$-hypergraph}
   & \parbox{0.7in}{$\alpha_{02}(m,n)$\newline $\alpha_{02}^\ast(m,n)$}
   & \parbox{0.7in}{$\alpha_{12}(m,n)$\newline $\alpha_{12}^\ast(m,n)$}
   & .\;\,.\;\,.\;\,.\;\,.\;\,.\;\,
   & \parbox{0.7in}{$\alpha_{k2}(m,n)$\newline $\alpha_{k2}^\ast(m,n)$} \\
\hline
\parbox{1.6in}{unordered $\mathfrak{q}$-hypergraph\newline without multiple edges}
   & \parbox{0.7in}{$\alpha_{03}(m,n)$\newline $\alpha_{03}^\ast(m,n)$}
   & \parbox{0.7in}{$\alpha_{13}(m,n)$\newline $\alpha_{13}^\ast(m,n)$}
   & .\;\,.\;\,.\;\,.\;\,.\;\,.\;\,
   & \parbox{0.7in}{$\alpha_{k3}(m,n)$\newline $\alpha_{k3}^\ast(m,n)$} \\
\hline
\parbox{1.6in}{unordered $\mathfrak{q}$-hypergraph}
   & \parbox{0.7in}{$\alpha_{04}(m,n)$\newline $\alpha_{04}^\ast(m,n)$}
   & \parbox{0.7in}{$\alpha_{14}(m,n)$\newline $\alpha_{14}^\ast(m,n)$}
   & .\;\,.\;\,.\;\,.\;\,.\;\,.\;\,
   & \parbox{0.7in}{$\alpha_{k4}(m,n)$\newline $\alpha_{k4}^\ast(m,n)$} \\
\hline
\end{tabular}
\end{center}}
\noindent
In the table
$\alpha_{ij}(m,n)$ [$\alpha_{ij}^\ast(m,n)$]
is the number of all hypergraphs [$T_0$-hyper\-graphs] that satisfy the property defined in the $j$-th row of the first
column and the property
$\mathfrak{q}_i$.
Here we suppose that properties
$\mathfrak{q}$ and $\mathfrak{q}_i$, $i\in\overline{0,k}$,
have sense for all hypergraphs, both in ordered and in unordered case. In the paper the tables will most often have
four columns, that is,
$k=3$,
and the properties
$\mathfrak{q}_i$, $i\in\overline{0,3}$,
will be, respectively, the properties `to be an arbitrary', `to be without empty edges', `to be without full edges',
and `to be without empty and full edges'. The table of the above-given type serves only for introducing new notations,
and as the first column is completely determined if the property
$\mathfrak{q}$
is given, then the table can be given in a more simplified form:
{\small\begin{center}
\begin{tabular}{|l|l|l|l|l|} \hline
            & $\mathfrak{q}_0$  & $\mathfrak{q}_1$  & .\;\,.\;\,.\;\,.\;\,.\;\,.\;\, & $\mathfrak{q}_k$ \\
\hline
to be $\mathfrak{q}$-hypergraph
            & $\alpha_{01}(m,n)$ & $\alpha_{11}(m,n)$ & .\;\,.\;\,.\;\,.\;\,.\;\,.\;\, & $\alpha_{k1}(m,n)$ \\
\hline
\end{tabular}
\end{center}}
\noindent
These tables (both a complete table and its simplified form) are called $\mathfrak{q}$-tables. If
$\alpha_{ij}(m,n)$
is a number in the given $\mathfrak{q}$-table, then we denote the corresponding class of hypergraphs by
$\mathfrak{H}(\alpha_{ij})$,
and the property (the combination of the corresponding properties) which determines the class
$\mathfrak{H}(\alpha_{ij})$ by $\mathfrak{p}(\alpha_{ij})$.

Let us agree that in the case when
$\mathfrak{q}$
or some of
$\mathfrak{q}_i$
is the property `to be arbitrary', the corresponding cell would be left empty. If the properties
$\mathfrak{q}_i$
are the same for several tables, then we shall ``glue'' these separate tables into a bigger one.
If the appropriate formulas for some of the numbers in the table are not adduced, and they are not
direct consequences of the given ones, the corresponding problems are not solved yet.

We usually get the number of unordered [ordered] $\mathfrak{p}$-hyper\-graphs without multiple edges
from the number of corresponding ordered [unordered] $\mathfrak{p}$-hyper\-graphs without multiple edges
by multiplying the latter number by
$1/m!$ [$m!$];
in that case we say that the {\em $G_2$-transformation\/} [{\em $G^{-1}_2$-transformation\/}] can be
applied on
$\mathfrak{p}$.

If the numbers
$\alpha_{ij}(m,n)$, $i\in\overline{0,n}$, $j\in\overline{2,4}$,
can be obtained from
$\alpha_{i1}(k,l)$
with the help of transformations
$G_1$, $G_2$ and $G_3=G_3'\,\circ\,G_2$
(the transformation
$G_3$
should be understood as consecutive application first of the transformation
$G_2$,
and then of the transformation
$G_3'$),
then we say that $\mathfrak{q}$-table is {\it regular\/}. It is clear that if $\mathfrak{q}$-table
is regular, then $T_0$-$\mathfrak{q}$-table is also regular. If the properties
$\mathfrak{q}$ and $\mathfrak{q}_i$, $i\in\overline{0,k}$,
allow $T_0$-filtration, then we say that the corresponding $\mathfrak{q}$-table is completely regular.
If $\mathfrak{q}$-table is completely regular, then the transformation
$F_0$
with each of the transformations
$G_i$, $i\in\overline{3}$,
commutates, i.e.,
$F_0\,\circ\,G_i=G_i\,\circ\,F_0$
for every
$i\in\overline{3}$,
on every combination of properties which is determined by the given $\mathfrak{q}$-table.

As we have already said, our task is to find numbers
$\alpha_{ij}^\ast(m,n)$.
Let us suppose that the properties
$\mathfrak{q}$ and $\mathfrak{q}_i$
allow $T_0$-filtration. Now if we know
$\alpha_{ij}(m,n)$,
then in order to calculate numbers
$\alpha_{ij}^\ast(m,n)$
it is sufficient just to use Theorem \ref{trm4}. So, in that case the problem of finding
$\alpha_{ij}^\ast(m,n)$
is considered solved if we know
$\alpha_{ij}(m,n)$.
In case when a $\mathfrak{q}$-table is completely regular, it follows from the given above considerations
that it is sufficient to know the numbers
$\alpha_{i1}(m,n)$,
and if we calculate them, the problem is considered solved.

Let us note that
$\alpha_{\mathfrak{p}}^\ast(0,n)=0$
for every
$n\ge2$, and $\alpha_{\mathfrak{p}}^\ast(0,1)=\alpha_{\mathfrak{p}}(0,1)$.
Hence, in the following we often assume that
$m\ge1$.

Let us illustrate all we have stated here with the help of a simple case.
Introduce the notations for the following classes of hypergraphs (here `wo.' stands for `without'):

\def\un#1{\underline{#1}}

{\small\begin{center}
\begin{tabular}{|c|c|c|c|c|} \hline
            & \parbox{0.9in}{\ } & \parbox{0.9in}{\scriptsize wo.\ empty edges}
            & \parbox{0.9in}{\scriptsize wo.\ full edges}
            & \parbox{0.9in}{\scriptsize wo.\ empty edges\newline and full edges} \\
\hline
            & $\alpha_{01}(m,n)$           & $\alpha_{11}(m,n)$
            & $\alpha_{21}(m,n)$           & $\alpha_{31}(m,n)$ \\\hline
no $\cap$-property
            & $\bar\alpha_{01}(m,n)$       & $\bar\alpha_{11}(m,n)$
            & $\bar\alpha_{21}(m,n)$       & $\bar\alpha_{31}(m,n)$ \\\hline
\end{tabular}
\end{center}}

\smallskip

\noindent
Also, introduce the following notations:
\begin{equation}
  \lambda_1(i,j)=[i]_j,\quad
  \lambda_2(i,j)=i^j,\quad
  \lambda_3(i,j)=C^j_i,\quad
  \lambda_4(i,j)=C^{j}_{i+j-1}\label{f3}
\end{equation}
for every
$i,j\in\text{\bf N}$;
here
$[i]_j$
is the falling factorial (see \cite{Agn}). By using Theorem \ref{trm4}, Theorem \ref{41}, and elementary
combinatorics we immediately get the following proposition:

\def\imtch#1#2#3#4{\hbox{\indent
                   \vtop{\hsize=#3cm\noindent #1}\hfill
                   \vtop{\hsize=.6cm\noindent --- }\vtop{\hsize=#4cm\noindent #2}\hfill}\smallskip}

\begin{proposition}\label{pp34}
  For every
  $j\in\overline{0,3}$ and $k\in\overline{4}$,
  $$
  \alpha_{jk}(m,n)=\lambda_k(2^n-\left[{\scriptstyle\frac{j+1}{2}}\right],m)\quad\text{and}\quad
  \alpha_{jk}^\ast(m,n)=\sum_{i=1}^{n}s_{n,i}\,\lambda_k(2^i-\left[{\scriptstyle\frac{j+1}{2}}\right],m),
  $$
  where
  $[x]$
  is the floor function.
\end{proposition}

Now, for example, as, obviously,
$\alpha_{02}^\ast(m,n)=[2^{m}]_n$,
then from Proposition \ref{pp34} we get the well known equality
$\sum_{i=1}^{n}(2^i)^ms_{n,i}=[2^{m}]_n$.

As {`no $\cap$-property'}-table is regular, it is sufficient to see that the following statement holds.

\begin{proposition}\label{pp35}
  For the numbers
  $\bar{\alpha}_{i1}(m,n)$, $i\in\overline{0,3}$,
  and for every
  $n\ge1$,
  we have that
  $$
  \bar{\alpha}_{i1}(m,n)=\alpha_{i1}(m,n)+\sum_{j=1}^{n-1}(-1)^jC^j_n\alpha_{2[i/2],1}(m,n-j)
  $$
  if
  $m>1$,
  and
  $\bar{\alpha}_{01}(1,n)=\bar{\alpha}_{21}(1,n)=1$
  and
  $\bar{\alpha}_{11}(1,n)=\bar{\alpha}_{31}(1,n)=0$.
\end{proposition}

Note that the functions
$\bar{\alpha}_{11}(m,n)$ and $\bar{\alpha}_{21}(m,n)$
are symmetric.

Denote by
$A_{1i}(x,y)$, $i\in\overline{4}$,
the exponential generating function for
$\alpha_{1i}(m,n)$.
Then by Theorem \ref{trm8} and Theorem \ref{trm9}, using Lah and Stirling transforms \cite{Slo2},
(see also \cite{Wlf}), it is not difficult to show that
{\small
$$
\begin{array}{l}
A_{11}(x,y)=\displaystyle\sum_{m\geq 0}\sum_{n\geq 0}\alpha _{11}(m,n)\,\frac{y^{m}}{m!}\,
            \frac{x^{n}}{n!}=\sum_{m\geq 0}\sum_{n\geq 0}[2^{n}-1]_{m}\,
            \frac{y^{m}}{m!}\,\frac{x^{n}}{n!}={}\\
\phantom{A_{11}(x,y)}=\displaystyle\sum_{n\geq 0}(1+y)^{2^{n}-1}\frac{x^{n}}{n!},\\
A_{12}(x,y)=\displaystyle\sum_{m\geq 0}\sum_{n\geq 0}\alpha _{12}(m,n)\,\frac{y^{m}}{m!}\
            \frac{x^{n}}{n!}=A_{11}(x,e^{y}-1)=\ \sum_{n\geq 0}e^{(2^{n}-1)y}\frac{x^{n}}{n!},\\
A_{13}(x,y)=\displaystyle\sum_{m\geq 0}\sum_{n\geq 0}\alpha _{13}(m,n)\ y^{m}\,\frac{x^{n}}{n!}=A_{11}(x,y)=
            \frac{1}{1+y}\cdot \sum_{n\geq 0}e^{2^{n}x}\,\frac{(\ln (1+y))^{n}}{n!},\\
A_{14}(x,y)=\displaystyle\sum_{m\geq 0}\sum_{n\geq 0}\alpha _{14}(m,n)\ y^{m}\,\frac{x^{n}}{n!}=A_{11}(x,\frac{y}{1-y})=
            \sum_{n\geq 0}(1-y)^{-2^{n}+1}\,\frac{x^{n}}{n!}=\\
\phantom{A_{14}(x,y)}=\displaystyle(1-y)\cdot \sum_{n\geq 0}e^{2^{n}x}\,\frac{(-\ln (1-y))^{n}}{n!}.
\end{array}
$$}

\section{Covers, $k_{\le}$-dimensional and $k$-uniform hypergraphs}

A cover $H$ is a {\em $k_{\le}$-cover\/} [$k$-{\em cover\/}],
$k\in\text{\bf N}$,
if for every vertex there are no more than $k$ [exactly $k$] edges containing it. Note again that the property
`to be without singular point' is equivalent to the property `to be cover with no $\cap$-property'. The dual property
of the property `to be an ordered $k_{\le}$-cover [$k$-cover]' is `to be a $k_{\le}$-dimensional hypergraph without
empty edges [a $k$-uniform hypergraph]'. Introduce the following notation:

\def\bbt{\vphantom{\bar{\bar{\bar{\beta}}}}\bar{\beta}}
\def\bbbt{\vphantom{\bar{\bar{\bar{\beta}}}}\bar{\bar{\beta}}}
\def\bbd{\underline{\beta}\,}
\def\bbbd{\underline{\underline{\beta}}\,}
\def\mmt{\vphantom{\bar{\bar{\bar{\mu}}}}\bar{\mu}}
\def\mmmt{\vphantom{\bar{\bar{\bar{\mu}}}}\bar{\bar{\mu}}}

{\small\begin{center}
\begin{tabular}{|l|l|l|l|l|} \hline
   & \parbox{0.8in}{\ } & \parbox{0.9in}{\scriptsize wo.\ empty edges}
   & \parbox{0.8in}{\scriptsize wo.\ full edges}
   & \parbox{0.9in}{\scriptsize wo.\ empty edges\newline and full edges} \\
\hline
cover
   & $\beta_{01}(m,n)$   & $\beta_{11}(m,n)$   & $\beta_{21}(m,n)$   & $\beta_{31}(m,n)$   \\\hline
$k$-cover
   & $\bbt_{01}(m,n,k)$  & $\bbt_{11}(m,n,k)$  & $\bbt_{21}(m,n,k)$  & $\bbt_{31}(m,n,k)$  \\\hline
$k_{\le}$-cover
   & $\bbbt_{01}(m,n,k)$ & $\bbbt_{11}(m,n,k)$ & $\bbbt_{21}(m,n,k)$ & $\bbbt_{31}(m,n,k)$ \\\hline
\parbox{1.2in}{\baselineskip9pt without singular\newline vertices}
   & $\beta_{41}(m,n)$   & $\beta_{51}(m,n)$   & $\beta_{61}(m,n)$   & $\beta_{71}(m,n)$   \\\hline
\end{tabular}
\end{center}}

A cover [$k_{\le}$-cover]
$H=(V,{\mathcal E})$, ${\mathcal E}=(e_\lambda,\lambda\in\Lambda)$,
is a {\em minimal cover\/} [{\em minimal $k_\le$-cover\/}] if for every
$\lambda\in\Lambda$
the subhypergraph
$H[V,\Lambda\backslash\{\lambda\}]$
is not a cover. A minimal cover [minimal $k_\le$-cover] $H$ never contains an empty edge and multiple edges,
and if
$|{\mathcal E}H|>1$,
it never contains a full edge either (if
$|{\mathcal E}H|=1$,
a cover contains only one full edge, and it is minimal). Denote by
$\mu_{01}(m,n)$ [$\mmt_{01}(m,n,k)$]
the number of all minimal covers [minimal $k_\le$-covers], and by
$\mu_{41}(m,n)$
the number of all minimal covers without $\cap$-property.

It is also possible to speak about minimal 1-covers, but it is easy to see that the number of such ordered
$(m,n)$-hypergraphs is
$S_{n,m}$,
and the number of ordered minimal $T_0$-1-covers having $m$ edges and $n$ vertices is
$m!$ if $m=n$, and $0$ if $m\ne n$.

\begin{theorem}\label{41}
  Each of the following properties allows $T_0$-filtration: `to be a cover', `to be a $k$-cover', `to be a
  $k_{\le}$-cover', `to be a minimal cover', `to be a minimal $k_\le$-cover'.
\end{theorem}

By using Theorems \ref{41}, \ref{trm4} and \ref{trm5} we have that all properties introduced by the above table and all the above
introduced minimal properties allow $T_0$-filtration. Let us recall that a cover without full edges is a proper
cover.

\begin{figure}[ht]
\begin{center}
\includegraphics{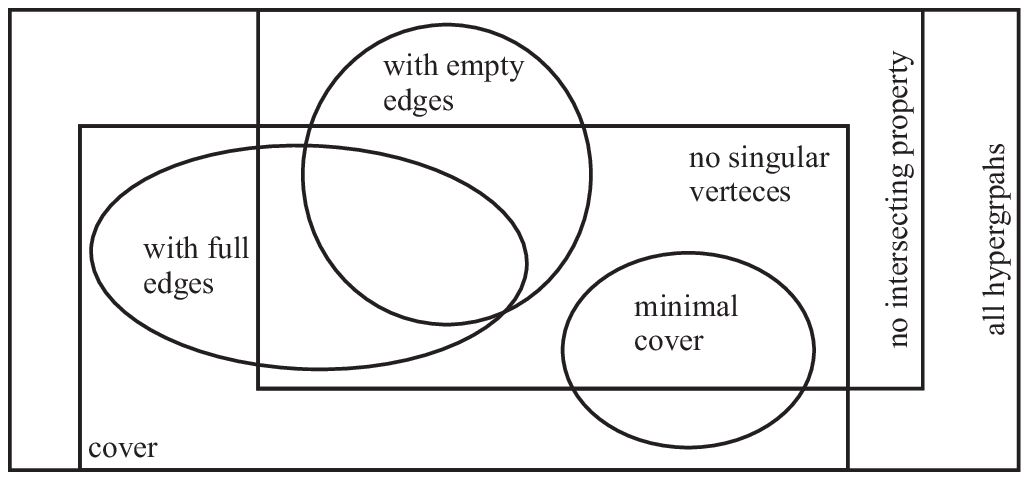}
\vspace{-0.2cm}
\end{center}
\caption{}
\label{fig:1}
\end{figure}

The relations between some of the introduced classes of ordered classes, when
$m\ge2$,
are given by Fig.\ \ref{fig:1}. The numbers introduced by the table are in some sense basic. For example,
if we know them, we are able to enumerate every constituent of the classes corresponding to these numbers.

Note that
$\bar{\bar\beta}_{ij}(m,n,1)=\bar\beta_{ij}(m,n,1)$
for every
$i\in\overline{0,3}$ and $j\in\overline{4}$, $\bar\beta_{01}(m,n,1)=\bar\beta_{02}(m,n,1)$,
$\bar\beta_{03}(m,n,1)=\bar\beta_{04}(m,n,1)$, $\bar\beta_{11}(m,n,1)=\bar\beta_{31}(m,n,1)=S_{n,m}$,
and, also,
$\beta_{01}(m,n)=(2^m-1)^n$.
Also note that `wo.\ singular vertices'-table is completely standard, and holds

\begin{proposition}\label{pp41}
  For the numbers
  $\beta_{i1}(m,n)$, $i\in\overline{4,7}$,
  and for every
  $n\ge1$,
  we have that
  $$
  \beta_{i1}(m,n)=\bar{\alpha}_{i1}(m,n)+\sum_{j=1}^{n-1}(-1)^jC^j_n\bar{\alpha}_{\nu(i),1}(m,n-j)
  $$
  if
  $m>1$,
  and
  $\beta_{i1}(1,n)=0$;
  here
  $\nu(i)=0$
  if $i$ is an even integer, and
  $\nu(i)=1$
  if $i$ is an odd one.
\end{proposition}

Note that the function
$\beta_{31}$
is symmetric. Also note that simpler formulas are possible, for example,
$$
\beta_{41}(m,n)=\sum_{i=0}^n(-1)^iC^i_n2^i[n-i]_m.
$$

Let
$H=(V,{\mathcal E})$
be a hypergraph, and let
$V_0(H)$
be the set of all isolated vertices of $H$. Denote by
$(H)_0$
the induced subhypergraph
$H[V\backslash V_0(H)]$.
Let
$\mathfrak{p}$
be a hypergraph property. We say that the property
$\mathfrak{p}$
is $\emptyset$-stable if for every hypergraph $H$ it holds that
$H\in\mathfrak{H}_{\mathfrak{p}}$
iff
$(H)_0\in\mathfrak{H}_{\mathfrak{p}}$,
i.e., iff
$\mathfrak{p}$
is invariant with respect to the operations of adding and cancelling isolated vertices.

\begin{theorem}\label{42}
  Let
  $\mathfrak{p}$
  be a\/ $\emptyset$-stable property that allows $T_0$-filtration, and let
  $\alpha_{\mathfrak{p}}(m,n)$
  be the number of all
  $(m,n)$-hypergraphs
  having property
  $\mathfrak{p}$.
  Then the number of all $T_0$-covers having property
  $\mathfrak{p}$
  is
  \begin{equation}
    \beta^\ast_{\mathfrak{p}}(m,n)=\sum^{n+1}_{i=1}s_{n+1,i}\,\alpha_{\mathfrak{p}}(m,i-1).\label{f4}
  \end{equation}
\end{theorem}

{\it Proof.\/}\;
Denote by
$\mathfrak{p}'$
the following hypergraph property: a hypergraph $H$ belongs to
$\mathfrak{H}_{\mathfrak{p}'}$ if $H\in\mathfrak{H}_{\mathfrak{p}}$ and $VH[1]$
is an isolated vertex. The property
$\mathfrak{p}'$
allows $T_0$-filtration. Denote by
$\gamma_{\mathfrak{p}}(m,n)$
the number of all hypergraphs from the set
$\mathfrak{H}_{\mathfrak{p}'}(m,n)$.
It is clear that
$\gamma_{\mathfrak{p}}(m,n)=\alpha_{\mathfrak{p}}(m,n-1)$.
By using Theorem \ref{trm4} we get
\begin{equation}
  \gamma^\ast_{\mathfrak{p}}(m,n)=\sum^n_{i=1}s_{n,i}\,\gamma_{\mathfrak{p}}(m,i)=
  \sum^n_{i=1}s_{n,i}\,\alpha_{\mathfrak{p}}(m,i-1).\label{f5}
\end{equation}
As a $T_0$-hypergraph has one isolated vertex at most, we have that
\begin{equation}
  \gamma^\ast_{\mathfrak{p}}(m,n)=
  \beta^\ast_{\mathfrak{p}}(m,n-1)\qquad\text{for every}\quad n\ge2.\label{f6}
\end{equation}
Now formula (\ref{f4}) follows from (\ref{f5}) and (\ref{f6}).\,$\square$

\smallskip
It is easy to see that the properties `to be a hypergraph' and `to be a hypergraph without empty edges' are
$\emptyset$-stable, and the properties `to be a hypergraph without full edges' and `to be a hypergraph without
empty and full edges' are not $\emptyset$-stable. Therefore, by using Theorem \ref{42} we get the first two formulas
of the next proposition.

\begin{proposition}\label{pp42}
  For the numbers
  $\beta_{j1}^\ast(m,n)$, $j\in\overline{0,3}$,
  we get the following formulas:
  $$
  \begin{array}{l}
     \beta_{01}^\ast(m,n)=\displaystyle\sum^{n+1}_{i=1}[2^{i-1}]_m\,s_{n+1,i},\quad
     \beta_{21}^\ast(m,n)=\beta_{01}^\ast(m,n)-m\alpha_{21}^\ast(m-1,n),\\
     \beta_{11}^\ast(m,n)=\displaystyle\sum^{n+1}_{i=1}[2^{i-1}-1]_m\,s_{n+1,i},\quad
     \beta_{31}^\ast(m,n)=\beta_{11}^\ast(m,n)-m\alpha_{31}^\ast(m-1,n).
  \end{array}
  $$
\end{proposition}

Note that the functions
$\beta_{11}^\ast(m,n)$ and $\beta_{21}^\ast(m,n)$
are symmetric. Also note that
$$
\beta^\ast_{02}(m,n)=\sum^{n+1}_{i=1}2^{m(i-1)}s_{n+1,i}=[2^m-1]_n.
$$

Let us give some further relations between the introduced classes. Note that `minimal cover'-table
and `minimal cover without intersecting property'-table are regular.

\begin{proposition}\label{pp43}
  For minimal covers we get
  $$
  \mu_{01}(m,n)=\sum^{n}_{i=m}C^i_n\,S_{i,m}\,m!\,(2^m-m-1)^{n-i}
  $$
  if $n\ge m$, and $\mu_{01}(m,n)=0$ if $m>n$.
\end{proposition}

\smallskip
Using Theorem \ref{trm4} we can obtain the formula for the number
$\mu^\ast_{01}(m,n)$,
which can be transformed into the form
$$
\mu^\ast_{01}(m,n)=n!\,C^{n-m}_{2^m-m-1}.
$$

Without applying Theorem \ref{trm4} the latter formula can be proved in the following way. An ordered
$T_0$-$(m,n)$-hypergraph $H$ is a minimal cover iff there exists
$V'\subseteq HV$, $|V'|=m$,
such that
$H[V']$
is a 1-cover, and in
$H[VH\backslash V']$
every vertex is covered by at least two edges. So if an $n$-set is fixed, and we want to construct an
ordered minimal $T_0$-$(m,n)$-cover $H$ on it, we have
$m!\,C^m_n$
possibilities for
$H[V']$,
and
$[2^m-m-1]_{n-m}$
possibilities for
$H[VH\backslash V']$.
Therefore we get
$$
\mu^\ast_{01}(m,n)=m!\,C^m_n[2^m-m-1]_{n-m}=n!\,C^{n-m}_{2^m-m-1}.
$$

\begin{proposition}\label{pp44}
  It holds that
  $$
  \mu_{41}(m,n)=\sum_{i=0}^{n-1}(-1)^iC^i_n\,\mu_{01}(m,n-i)
  $$
  if $n\ge m$ and $m>1$, and $\mu_{41}(m,n)=0$ if either $n<m$ or $m=1$.
\end{proposition}

Note that for every
$k\ge n$,
$\bar\mu_{01}(m,n,k)=\mu_{01}(m,n)$ and $\bar\mu_{41}(m,n,k)=\mu_{41}(m,n)$.

Now consider the problem of finding the numbers
$\delta_{ij}(m,n,k)$,
where
$i\in\overline{0,3}$, $j\in\overline{1,4}$, and $\delta$
stands for $\bbt$, $\bbbt$, $\bbt^\ast$ or $\bbbt^\ast$.
Let us give an example, considering the problem of enumeration of all unordered
$T_0$-$(m,n)$-hypergraphs without multiple edges and without empty edges that are $2$-covers
[$2_\le$-covers], i.e., considering the problem of finding the number
$\bbt_{13}^\ast(m,n,2)$ (see \cite{Gou}) [$\bbbt_{13}^\ast(m,n,2)$].

\smallskip
{\bf Example 4.1.}\ Note that the dual property of the property `to be an ordered $T_0$-$2$-cover
without multiple edges' is the
property `to be an ordered $2$-uniform hypergraph without multiple edges, without any $2$-components, and
without any isolated vertex'. Let us find the number
$\bar{\theta}_{03}^\circ(m,n)$
of all unordered 2-uniform $(m,n)$-hypergraphs without multiple edges, and without 2-components, i.e., graphs
with $m$ edges and $n$ vertices, and without any 2-components.

The number of all graphs with $n$ vertices and $m$ edges is
$C^m_{C^2_n}$.
The maximum number of 2-components is, obviously,
$\min\{[n/2],m\}$.
The number of all graphs having the given $k$ $2$-components is
$C^{m-k}_{C^2_{n-2k}}$.
Disjoint $k$ $2$-components can be chosen in
$[C^{2k}_n(2k)!]/[(2!)^kk!]=[n]_{2k}/[(2!)^kk!]$
ways. Thus, by using the formula of inclusion and exclusion we get the formula
$$
\bar{\theta}_{03}^\circ(m,n)=
\sum_{k=0}^{\min\{[n/2],m\}}(-1)^k\,\dfrac{[n]_{2k}}{(2!)^kk!}\,C^{m-k}_{C^2_{n-2k}}.
$$

Similarly, for the number
$\bar{\bar{\theta}}_{03}^\circ(m,n)$
of all unordered $2_{\le}$-dimensional $(m,n)$-hypergraphs without multiple edges, and without 2-components,
i.e., graphs which have $m$ edges (among them there can be even loops) and $n$ vertices, and which do not
have any 2-components, we have the following formula:
$$
\bar{\bar{\theta}}_{03}^\circ(m,n)=
\sum_{k=0}^{\min\{[n/2],m\}}(-1)^k\,\dfrac{[n]_{2k}}{(2!)^kk!}\,C^{m-k}_{C^1_{n-2k}+C^2_{n-2k}}.
$$

The generative function for
$\bar{\theta}_{03}^\circ(m,n)$ and $\bar{\bar{\theta}}_{03}^\circ(m,n)$
can be found in \cite{Slo1}.

Denote by
$\bar{\theta}_{13}^\circ(m,n)$ [$\bar{\bar{\theta}}_{13}^\circ(m,n)$]
the number of all unordered 2-uniform [$2_{\le}$-dimen\-sional] $(m,n)$-covers without multiple edges [without
multiple edges, without empty edges], without 2-components, and without isolated vertices. Then we have
$$
\bar{\theta}_{13}^\circ(m,n)=\sum_{i=0}^{n}(-1)^iC^{i}_{n}\bar{\theta}_{03}^\circ(m,n-i),\quad
\bar{\bar{\theta}}_{13}^\circ(m,n)=\sum_{i=0}^{n}(-1)^iC^{i}_{n}\bar{\bar{\theta}}_{03}^\circ(m,n-i),
$$
and we get
$$
\bbt_{13}^\ast(m,n,2)=\dfrac{n!}{m!}\,\bar{\theta}_{13}^\circ(n,m),\quad
\bbbt_{13}^\ast(m,n,2)=\dfrac{n!}{m!}\,\bar{\bar{\theta}}_{13}^\circ(n,m).\,\square
$$

Let
$\delta$
replace here one of the letters
$\bbt$, $\bbbt$, $\bbt^\ast$ and $\bbbt^\ast$.
It is not difficult to see that for every
$j\in\overline{4}$
we can calculate all the numbers
$\delta_{ij}(m,n,k)$, $i\in\overline{3}$,
if we know the number
$\delta_{0j}(m,n,k)$.
We immediately obtain the number
$\delta_{03}(m,n,k)$
from the number
$\delta_{01}(m,n,k)$,
and the numbers
$\delta_{01}(m,n,k)$ and $\delta_{02}(m,n,k)$
can be obtained by considering the corresponding dual properties. Because of that we introduce
the following notation:

\def\thett{\vphantom{\bar{\bar{\bar{\theta}}}}\bar{\theta}}
\def\thttt{\vphantom{\bar{\bar{\bar{\theta}}}}\bar{\bar{\theta}}}

{\small\begin{center}
\begin{tabular}{|p{1.8in}|l|l|l|} \hline
          &                        & cover                & minimal cover\\\hline
$k$-uniform
          & $\theta_{01}(m,n,k)$   & $\theta_{11}(m,n,k)$ & $\theta_{21}(m,n,k)$  \\\hline
$k_{\le}$-dim.\ wo.\ empty edges
          & $\thett_{01}(m,n,k)$   & $\thett_{11}(m,n,k)$ & $\thett_{21}(m,n,k)$  \\\hline
$k$-uniform and $\neg\cap$-property
          & $\theta_{31}(m,n,k)$   & $\theta_{41}(m,n,k)$ & $\theta_{51}(m,n,k)$ \\\hline
{$k_{\le}$-dim.\ wo.\ empty edges\newline and $\neg\cap$-property}
          & $\thett_{31}(m,n,k)$   & $\thett_{51}(m,n,k)$ & $\thett_{51}(m,n,k)$ \\\hline
\end{tabular}
\end{center}}

\begin{figure}
\begin{center}
\includegraphics{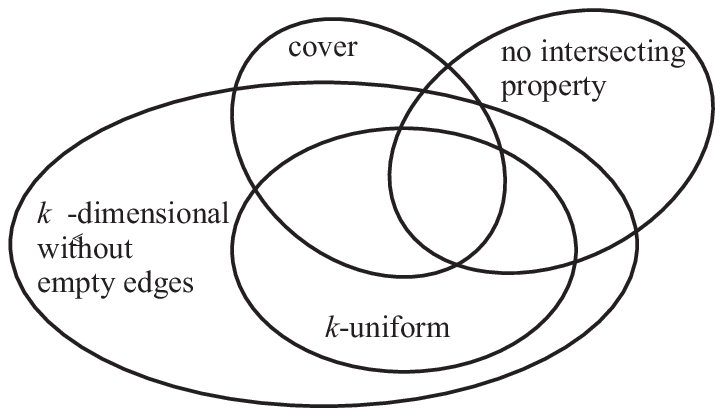}
\vspace{-0.2cm}
\end{center}
\caption{}
\label{fig:2}
\end{figure}

Note that the property `to be a $k_{\le}$-dimensional hypergraph' does not satisfy the conditions
from Theorem \ref{trm4} and Theorem \ref{trm5}. Really, the formula from Theorem \ref{trm4} does not
give the number of all $T_0$-$k_\le$-hypergraphs. This also holds for the class of all $k$-dimensional
hypergraphs and the class of all $k$-uniform hypergraphs. But it is not difficult to see that
$$
\theta_{01}(m,n,k)=[C^k_n]_m,\quad\theta_{11}(m,n,k)=\sum_{i=0}^{n-k}(-1)^iC^i_n[C^k_{n-i}]_m,
$$
and for every
$j\in\overline{0,2}$,
$$
\theta_{3+j,1}(m,n,k)=\sum_{i=0}^{k-1}(-1)^{i}C^i_n\theta_{j1}(m,n-i,k-i)
$$
if $m>1$, and $\theta_{3+j,1}(1,n,k)=0$. Also,
$$
\bar\theta_{01}(m,n,k)=\displaystyle[\bar{C}^k_n]_m,\quad
\bar\theta_{11}(m,n,k)=\sum_{i=0}^{n-1}(-1)^iC^i_n[\bar{C}^k_{n-i}]_m;
$$
here
$\bar C^{s}_{t}=\displaystyle\sum_{i=1}^sC^i_t$.
Introduce the numbers
$$
\theta_{01}'(m,n,k)=[C^k_n+1]_m,\quad\theta_{11}'(m,n,k)=\sum_{i=0}^{n-k}(-1)^iC^i_n[C^k_{n-i}+1]_m,
$$
Then we get for every
$j\in\overline{0,1}$,
$$
\bar\theta_{3+j,1}(m,n,k)=\bar\theta_{j1}(m,n,k)+\sum_{i=1}^{k-1}(-1)^{i}C^i_n\bar\theta_{j1}'(m,n-i,k-i)
$$
if $m>1$, and $\theta_{3+j,1}(1,n,k)=0$.
Also we get
$$
\bar\theta_{51}(m,n,k)=\sum_{i=0}^{k-1}(-1)^{i}C^i_n\bar\theta_{21}(m,n-i,k-i).
$$
The problem of finding the number
$\theta_{21}(m,n,k)$ (and $\theta_{51}(m,n,k)$)
is unsolved yet. But we have the following formulas:
$$
\delta^\ast_{ij}(m,n+1,k)=\delta^\ast_{i-1,j}(m,n+1,k)-(n+1)\,\delta^\ast_{1i}(m,n,k),\quad
i=1,4;\;j\in\overline{4},
$$
where
$\delta$
stands here for
$\theta$ or $\bar\theta$,
with help of which we pass from the numbers of the second column to the numbers of the first
column of the given table. From the numbers of the third column to the numbers of the second
column of the given table we pass with help of the following formulas:
$$
\theta^\ast_{i1}(m,n,k)=[n]_m\theta^\ast_{i-1,2}(m,n-m,k-1),\quad i=2,5,
$$
and
$$
\bar\theta^\ast_{i1}(m,n,k)=[n]_m\sum_{j=1}^{m}C^j_m\,\bar\theta^\ast_{i-1,2}(j,n-m,k-1),\quad i=2,5
$$
(it is clear that
$m\le n$).
Finally, we have the formulas
$$
\delta^\ast_{3i}(m,n+1,k)=\delta^\ast_{0i}(m,n+1,k)-(n+1)\,\delta^\ast_{3i}(m,n,k),\quad
i\in\overline{4};
$$
here
$\delta$
stands for
$\theta$ or $\bar\theta$.
Therefore, the problem of finding the numbers
$\theta^\ast_{ij}(m,n,k)$  and $\bar\theta^\ast_{ij}(m,n,k)$,
that are defined by the above given table, can be reduced to the same problem for the numbers
$\theta^\ast_{0s}(m,n,k)$, $s\in\overline{4}$, and $\bar\theta^\ast_{0s}(m,n,k)$, $s\in\overline{4}$,
and in order to find the latter numbers we have to make use of the formula from Theorem \ref{trm2}.

Let
$\alpha=(\alpha_1,\dots,\alpha_n)$ and $\beta=(\beta_1,\dots,\beta_n)$
be arbitrary $n$-tuples.
We write
$\beta\le\alpha$, if $\beta_i\le\alpha_i$
for every
$i\in\overline{n}$.

Take a partition
$\pi=\{V_1,\dots,V_i\}\in\Pi(n)$
of the $n$-set $W_n$, and fix a $k$,
$k<n$.
Denote by
$\nu(\pi,k)$
the number of all $k$-sets
$V\subseteq W_n$
which can be represented in the form
$V=V_{t_1}\cup\dots\cup V_{t_j}$
for some
$1\le t_1<\dots<t_j\le i$.
Let
$\alpha=(\alpha_1,\dots,\alpha_n)$
be a partition type, and let
$\pi\in\Pi(\alpha_1,\dots,\alpha_n)$.
Then it is easy to get
$$
\nu(\pi,k)=\nu(\alpha,k)=\sum_{\beta=(\beta_1,\dots,\beta_n),\,\sigma(\beta)=k,\,\beta\le\alpha}
C_{\alpha_1}^{\beta_1}\cdot\dots\cdot C_{\alpha_n}^{\beta_n}.
$$
Now by using Theorem \ref{trm2} we obtain the following theorem.

\begin{theorem}\label{43}
  For every
  $i\in\overline{4}$,
  $$
  \theta^\ast_{0i}(m,n,k)=
  \sum_{\alpha=(\alpha_1,\dots,\alpha_n),\,\sigma(\alpha)=n}(-1)^{n-\alpha_1-\dots-\alpha_n}
  c(\alpha)\,\lambda_i(\nu(\alpha,k),m).
  $$
\end{theorem}

Also, if we introduce the number
$$
\nu_\le(\alpha,k)=\sum_{\beta=(\beta_1,\dots,\beta_n),\,0<\sigma(\beta)\le k,\,\beta\le\alpha}
C_{\alpha_1}^{\beta_1}\cdot\dots\cdot C_{\alpha_n}^{\beta_n},
$$
we get

\begin{theorem}\label{44}
  For every
  $i\in\overline{4}$,
  $$
  \bar\theta^\ast_{0i}(m,n,k)=\sum_{\alpha=(\alpha_1,\dots,\alpha_n),\,\sigma(\alpha)=n}(-1)^{n-\alpha_1-\dots-\alpha_n}
  c(\alpha)\,\lambda_i(\nu_\le(\alpha,k),m).
  $$
\end{theorem}

Note that any of the transformations
$G_i$, $i\in\overline{3}$,
can be applied on
$\mathfrak{p}(\theta^\ast_{01})$ and $\mathfrak{p}(\bar\theta^\ast_{01})$,
and by using them the corresponding formulas for the numbers
$\theta^\ast_{0i}$ and $\bar\theta^\ast_{0i}$, $i\in\overline{2,4}$
can be obtained from the above formulas for the numbers
$\theta^\ast_{01}$ and $\bar\theta^\ast_{01}$.
So obtained formulas would be equivalent to the formulas given by the above two theorems.

\section{Connected hypergraphs}

Let us begin with the following theorem:

\begin{theorem}\label{51}
  The property `to be connected' allows $T_0$-filtration.
\end{theorem}

It is easy to see that the property `to be a connected $k$-uniform [$k_\le$-uniform] hypergraph'
does not allow $T_0$-filtration. Also note that any connected hypergraph is a cover. Obviously,
the opposite does not hold. However, any cover with intersecting property is connected. Introduce
the notations:\par\smallskip

{\small\begin{center}
\begin{tabular}{|l|l|l|l|l|} \hline
& & \parbox{0.9in}{\scriptsize wo.\ empty edges}
  & \parbox{0.9in}{\scriptsize wo.\ full edges}
  & \parbox{0.9in}{\scriptsize wo.\ full edges and wo.\ empty edges} \\
\hline
\parbox{1.1in}{\scriptsize connected}
  &  $\omega_{01}(m,n)$ & $\omega_{11}(m,n)$ & $\omega_{21}(m,n)$ & $\omega_{31}(m,n)$ \\\hline
\parbox{1.1in}{\scriptsize connected $k$-uniform}
  &  $\bar\omega_{01}(m,n,k)$ & & & \\\hline
\parbox{1.1in}{\scriptsize connected $k_\le$-dim.}
  &  $\bar{\bar\omega}_{01}(m,n,k)$ & $\bar{\bar\omega}_{11}(m,n,k)$ & & \\\hline\hline
\parbox{1.1in}{\scriptsize connected without\newline $\cap$\,-\,property}
  &  $\omega_{41}(m,n)$ & $\omega_{51}(m,n)$ & $\omega_{61}(m,n)$ & $\omega_{71}(m,n)$ \\\hline
\parbox{1.1in}{\scriptsize connected $k$-uniform\newline without $\cap$\,-\,property}
  &  $\bar\omega_{41}(m,n,k)$ & & & \\\hline
\parbox{1.1in}{\scriptsize connected $k_\le$-dim.\newline without $\cap$\,-\,property}
  &  $\bar{\bar\omega}_{41}(m,n,k)$ & $\bar{\bar\omega}_{51}(m,n,k)$ & & \\\hline
\end{tabular}
\end{center}}

The numbers
$\omega_{1j}(m,n)$, $j\in\overline{4}$,
are calculated in Proposition \ref{pp57}. The next four propositions show how from them we
can obtain the numbers
$\omega_{ij}(m,n)$, $i\in\overline{0,7}$, $i\ne1$, $j\in\overline{4}$.

\begin{proposition}\label{pp51}
  For the numbers
  $\omega_{0k}(m,n)$, $k\in\overline{4}$,
  we get\newline
  $
  \begin{array}{l}
  \hskip1cm\omega_{01}(m,n)=m\,\omega_{11}(m-1,n)+\omega_{11}(m,n), \\
  \hskip1cm\omega_{02}(m,n)=\sum_{i=0}^{m-1}C^{i}_{m}\,\omega_{12}(m-i,n),\\
  \hskip1cm\omega_{03}(m,n)=\omega_{13}(m-1,n)+\omega_{13}(m,n), \\
  \hskip1cm\omega_{04}(m,n)=\sum_{i=0}^{m-1}\omega_{14}(m-i,n);
  \end{array}
  $\newline
  suppose that for every
  $n>1$, $\omega_{0k}(0,n)=0$ and $\omega_{0k}(1,n)=1$,
  and, also,
  $\omega_{0k}(0,1)=1$ and $\omega_{0k}(1,1)=2$.
\end{proposition}

\begin{proposition}\label{pp52}
  For the numbers
  $\omega_{jk}(m,n)$, $j\in\overline{2,3}$, $k\in\overline{4}$,
  we get\newline
  $
  \begin{array}{l}
  \hskip1cm\omega_{j1}(m,n)=\omega_{j-2,1}(m,n)-m\alpha_{j1}(m-1,n), \\
  \hskip1cm\omega_{j2}(m,n)=\omega_{j-2,2}(m,n)-\sum_{i=1}^mC^i_m\alpha_{j2}(m-i,n) \\
  \hskip1cm\omega_{j3}(m,n)=\omega_{j-2,3}(m,n)-\alpha_{j3}(m-1,n), \\
  \hskip1cm\omega_{j4}(m,n)=\omega_{j-2,4}(m,n)-\sum_{i=1}^{m}\alpha_{j4}(m-i,n);
  \end{array}
  $\newline
  here
  $\alpha_{j2}(0,n)=\alpha_{j4}(0,n)=1$.
\end{proposition}

\begin{proposition}\label{pp53}
  For the numbers
  $\omega_{jk}(m,n)$, $j\in\overline{4,5}$, $k\in\overline{4}$,
  we get
  $$
  \omega_{jk}(m,n)=\omega_{j-4,k}(m,n)-\sum_{i=1}^n(-1)^iC^i_n\beta_{0k}(m,n-i);
  $$
  here
  $\beta_{01}(m,0)=\beta_{03}(m,0)=0$ and $\beta_{02}(m,0)=\beta_{04}(m,0)=1$.
\end{proposition}

\begin{proposition}\label{pp54}
  For the numbers
  $\omega_{jk}(m,n)$, $j\in\overline{6,7}$, $k\in\overline{4}$,
  we get
  $$
  \omega_{jk}(m,n)=\omega_{j-4,k}(m,n)-\sum_{i=1}^n(-1)^iC^i_n\beta_{2k}(m,n-i);
  $$
  here
  $\beta_{21}(m,0)=\beta_{23}(m,0)=0$ and $\beta_{22}(m,0)=\beta_{24}(m,0)=1$.
\end{proposition}

And, of course, Theorem \ref{51} and Proposition \ref{pp33} imply

\begin{proposition}\label{pp55}
  For the numbers
  $\omega^\ast_{jk}(m,n)$, $j\in\overline{7}$, $k\in\overline{4}$,
  we get
  $$
  \omega^\ast_{jk}(m,n)=\sum_{i=1}^ns_{n,i}\omega_{jk}(m,i).
  $$
\end{proposition}

Let
$H=(V,{\mathcal E})$, ${\mathcal E}=(e_\lambda,\lambda\in\Lambda)$,
be a
$(m,n)$-hypergraph.
The bipartite graph
$K_{m,n}(H)=(V\cup\Lambda,E)$, $E=\{(v,\lambda)\in V\times\Lambda\,|\,v\in e_\lambda\}$,
with the partition classes $V$ and
$\Lambda$,
is called the {\it graph of incidence\/} of $H$. It is easy to see that an unordered [ordered] hypergraph $H$ is
connected iff the whole $n$-block of
$K_{m,n}(H)$
belongs to the same component of connectedness of this bipartite graph. Now the next proposition follows immediately.

\begin{proposition}\label{pp56}
  A hypergraph $H$ is an ordered $[$unordered\/$]$ connected $(m,n)$-hypergraph without empty edges iff
  $K_{m,n}(H)$
  is a connected graph.
\end{proposition}

\begin{corollary}\label{crll51}
  For all
  $m,n\in\text{\bf N}$, $\omega_{12}(m,n)=\omega_{12}(n,m)$,
  i.e.\ the function
  $\omega_{12}(m,n)$
  is symmetric.
\end{corollary}

Let
$\mathfrak{p}$
be a property of ordered [unordered] hypergraphs without empty edges. Let
$H=(V,{\mathcal E})$, ${\mathcal E}=(e_\lambda,\lambda\in\Lambda)$,
be an arbitrary ordered [unordered] hypergraph without empty edges. Suppose that there exist subsets
$V_1$ and $V_2$
of $V$ such that
$V_1\cup V_2=V$, $V_1\cap V_2=\emptyset$,
and
$e_\lambda\subseteq V_1$ or $e_\lambda\subseteq V_1$
for every
$\lambda\in\Lambda$,
and let
$H_1=H[V_1]$ and $H_2=H[V_2]$.
We say that the pair
$(H_1,H_2)$
is a {\em $\gamma$-decomposition\/} of $H$, and we write
$H=H_1\lor H_2$, if
$H_1$
is a connected hypergraph. We say that the property
$\mathfrak{p}$
is {\em invariant with respect to $\gamma$-decompositions\/} if for every
$H_1$ and $H_2$
such that
$H=H_1\lor H_2$,
it hold that
$H\in\mathfrak{H}_{\mathfrak{p}}$
iff
$H_1\in\mathfrak{H}_{\mathfrak{p}}$ and $H_2\in\mathfrak{H}_{\mathfrak{p}}$.

\begin{theorem}\label{52}
  Let
  $\mathfrak{p}$
  be a property that is invariant with respect to $\gamma$-decom\-positions, and let
  $\alpha_{\mathfrak{p}}(m,n)=|\mathfrak{H}_{\mathfrak{p}}(m,n)|$,
  where
  $n\ge2$.
  Denote by
  $\omega_{\mathfrak{p}}(m,n)$
  the number of all connected hypergraphs from
  $\mathfrak{H}_{\mathfrak{p}}(m,n)$,
  and by
  $\alpha_{\mathfrak{p}}'(m,n)$
  the number of all hypergraphs $H$ from
  $\mathfrak{H}_{\mathfrak{p}}(m,n)$
  for which the vertex
  $VH[1]$
  is an isolated vertex. Then
  $$
  \omega_{\mathfrak{p}}(m,n)=\alpha_{\mathfrak{p}}(m,n)-\alpha_{\mathfrak{p}}'(m,n)-
  \sum_{i=1}^m\sum_{j=1}^{n-1}\hat{\nu}_{\mathfrak{p}}(m,i)\,C^{j-1}_{n-1}\,
  \alpha_{\mathfrak{p}}(m-i,n-j)\,\omega_{\mathfrak{p}}(i,j);
  $$
  here
  $\hat{\nu}_{\mathfrak{p}}(m,i)=C^i_m$
  if $\mathfrak{p}$ is a property of ordered hypergraphs, and
  $\hat{\nu}_{\mathfrak{p}}(m,i)=1$
  if $\mathfrak{p}$ is a property of unordered hypergraphs.
\end{theorem}

{\it Proof.\/}\;
Denote by
$\hat\omega_{\mathfrak{p}}(m,n)$
the number of all disconnected hypergraphs from the class
$\mathfrak{H}_{\mathfrak{p}}(m,n)$.
Obviously,
\begin{equation}
  \omega_{\mathfrak{p}}(m,n)=\alpha_{\mathfrak{p}}(m,n)-\hat\omega_{\mathfrak{p}}(m,n).\label{f7}
\end{equation}
Let
$H=(V,{\mathcal E})$, ${\mathcal E}=(e_\lambda,\lambda\in\Lambda)$,
be a disconnected
$(m,n)$-hypergraph
from the class
$\mathfrak{H}_{\mathfrak{p}}(m,n)$. If $VH[1]$
is an isolated vertex, the hypergraph $H$ is a disconnected hypergraph, and the number of such hypergraphs is given
by the number
$\alpha_{\mathfrak{p}}'(m,n)$.
So, we can suppose that
$VH[1]$
is not an isolated vertex. Denote by
$V_0\subseteq V$
the set of vertices that belong to the component of connectedness containing the vertex
$v_1$,
and by
${\mathcal E}_0$
the subfamily of
${\mathcal E}$
which contains all edges
$e_\lambda$
such that
$e_\lambda\cap V_0\ne\emptyset$.
It is obvious that
$V_0\cap e=e$
for every
$e\in{\mathcal E}_0$,
and that
$H[V_0,\Lambda({\mathcal E}_0)]$
is a connected $(i,j)$-hypergraph from the class
$\mathfrak{H}_{\mathfrak{p}}(i,j)$,
where
$i=|{\mathcal E}_0|$ and $j=|V_0|$.
Also it is obvious that
$V_0\cap e=\emptyset$
for every
$e\in{\mathcal E}\backslash{\mathcal E}_0$.
Therefore we have that
\begin{equation}
  \hat\omega_{\mathfrak{p}}(m,n)=\alpha_{\mathfrak{p}}'(m,n)+
  \sum_{i=1}^m\sum_{j=1}^{n-1}\hat{\nu}_\mathfrak{p}(m,i)\,C^{j-1}_{n-1}\,
  \alpha_{\mathfrak{p}}(m-i,n-j)\,\omega_{\mathfrak{p}}(i,j).\label{f8}
\end{equation}
Now the formula follows from (\ref{f7}) and (\ref{f8}).\,$\square$

\smallskip
If a property
$\mathfrak{p}$
is invariant with respect to $\gamma$-decompositions, and, consequently, we can apply Theorem \ref{52},
we say that the {\em $F_1$-transformation\/} can be applied on
$\mathfrak{p}$.
Denote by
$F_1(\mathfrak{H}_{\mathfrak{p}})$
the class of all connected hypergraphs from
$\mathfrak{H}_{\mathfrak{p}}$.
Note that the property `to be a $T_0$-hypergraph' is invariant with respect to $\gamma$-decompositions.
Now if both the $F_0$-transformation and $F_1$-transformation can be applied on
$\mathfrak{p}$,
we have that they commutate on
$\mathfrak{p}$, i.e., $F_0\,\circ\,F_1(\mathfrak{H}_{\mathfrak{p}})=F_1\,\circ\,F_0(\mathfrak{H}_{\mathfrak{p}})$.
So we have that
$F_0[F_1(\mathfrak{H}(\omega_{1j})))]=F_1[F_0(\mathfrak{H}(\omega_{1j}))]$
for every
$j\in\overline{4}$.

Let us agree that below the symbol
$\nu_k(i,j)$
means
$C^j_i$ if $k=1,2$,
and
$\nu_k(i,j)=1$ if $k=3,4$.
Then from Theorem \ref{52} we obtain that

\begin{proposition}\label{pp57}
  For every
  $k\in\overline{4}$, $m\ge1$, and $n\ge2$,
  $$
  \omega_{1k}(m,n)=\lambda_k(2^{n}-1,m)-
  \sum_{i=0}^m\sum_{j=1}^{n-1}\nu_k(m,i)\,C^{j-1}_{n-1}\,\lambda_k(2^{n-j}-1,m-i)\,\omega_{1k}(i,j);
  $$
  here,
  $\omega_{1k}(0,1)=1$, and $\omega_{1k}(0,i)=0$ for all $i>1$.
\end{proposition}

Note that the numbers
$\omega_{12}(m,n)$
are considered, in a somewhat different context, in \cite{Kre}. Also let us note that on the basis of
a well-known (now proverbial) connection between exponential generative functions for graphs [hypergraphs]
and connected [hypergraphs] \cite{Har}, by using functions
$A_{1i}(x,y)$
we get
$$
\Omega _{1i}(x,y)=1+\ln(A_{1i}(x,y))
$$
for every $i\in\overline{4}$. So, for example, for
$i=2$
we get
$$
\begin{array}{l}
  \Omega_{12}(x,y)=\displaystyle1+
  \ln\left(\sum_{n\geq 0}e^{(2^{n}-1)y}\dfrac{x^{n}}{n!}\right)=1+e^yx+(e^{3y}-e^{2y})\dfrac{x^2}{2!}\,+\\
  \hskip1cm+\,(e^{7y}-3e^{4y}+2e^{3y})\dfrac{x^3}{3!}+(e^{15y}-4e^{8y}-3e^{6y}+12e^{5y}-6e^{4y})\dfrac{x^4}{4!}+\dots
\end{array}
$$
or
$$
\begin{array}{l}
\omega _{12}(1,n)=1,\quad \omega _{12}(2,n)=3^n-2^n,\quad \omega _{12}(3,n)=7^n-3\cdot4^n+2\cdot3^n,\\
\omega _{12}(4,n)=15^n-4\cdot8^n-3\cdot6^n+12\cdot5^n-6\cdot4^n,\quad\dots.
\end{array}
$$

Let
$\mathfrak{p}_1$ and $\mathfrak{p}_2$
be two properties that are invariant with respect to $\gamma$-decom\-positions. Then the property
$\mathfrak{p}_1\land\mathfrak{p}_2$
is also invariant with respect to $\gamma$-decompositions. It is obvious that the property
`to be a $T_0$-hypergraph' is a property that is invariant with respect to $\gamma$-decompositions.
Now we can find the number of all connected $k$-uniform $T_0$-hypergraphs and all connected
$k_{\le}$-dimensional $T_0$-hypergraphs.

\begin{proposition}\label{pp58}
  For every
  $m\ge1$, $n\ge2$, and $s\in\overline{4}$,
  $$
  \begin{array}{l}
    \displaystyle\bar{\omega}_{0s}^\ast(m,n,k)=\theta_{0s}^\ast(m,n,k)-\theta_{1s}^\ast(m,n-1,k)-\\
    \hskip3cm-\displaystyle
    \sum_{i=1}^m\sum_{j=1}^{n-1}\nu_k(m,i)\,C^{j-1}_{n-1}\,\theta_{0s}^\ast(m-i,n-j,k)\,\bar{\omega}_{0s}^\ast(i,j,k);
  \end{array}
  $$
  where
  $\bar{\omega}_{0s}^\ast(m,1,k)=1$
  if
  $k=1\land((m=1\land s=1,3)\lor(m\ge1\land s=2,4))$,
  and
  $\bar{\omega}_{0s}^\ast(m,1,k)=0$,
  otherwise;
  $\theta_{0s}^\ast(0,1,k)=1$
  and
  $\theta_{0s}^\ast(0,a,k)=1$ if $a>1$.
\end{proposition}

\begin{proposition}\label{pp59}
  For every
  $m\ge1$, $n\ge2$, and $s\in\overline{4}$,
  $$
  \begin{array}{l}
    \displaystyle\bar{\bar\omega}_{1s}^\ast(m,n,k)=\bar\theta_{0s}^\ast(m,n,k)-\bar\theta_{1s}^\ast(m,n-1,k)-\\
    \hskip3cm-\displaystyle\sum_{i=1}^m\sum_{j=1}^{n-1}\nu_k(m,i)\,C^{j-1}_{n-1}\,\bar\theta_{1s}^\ast(m-i,n-j,k)\,
    \bar{\bar\omega}_{0s}^\ast(i,j,k);
  \end{array}
  $$
  where
  $\bar{\bar\omega}_{1s}^\ast(m,1,k)=1$
  if
  $(m=1\land s=1,3)\lor(m\ge1\land s=2,4)$,
  and
  $\bar{\bar\omega}_{1s}^\ast(m,1,k)=0$,
  otherwise;
  $\bar\theta_{0s}^\ast(0,1,k)=1$
  and
  $\bar\theta_{0s}^\ast(0,a,k)=1$ if $a>1$.
\end{proposition}

In the end, let us remark that in addition to the properties discussed in the paper there are other
interesting properties which allow $T_0$-filtration. For example, each of the following properties
allows $T_0$-filtration: `to be an antichain' (i.e., `to be a dually $T_1$-hypergraph'), `to be a dually
$T_2$-hypergraph', `to be a hypergraph having $k$-intersecting property', etc. The case of $k$-intersecting
hypergraphs is considered in \cite{JK1}. The case of antichains is considered in \cite{KJ1}, and
the case of multiantichains is considered in \cite{KJ2}.


\begin{thebibliography}{}
\bibitem{Gou} I.\ P.\ Goulden  and D.\ M.\ Jackson, {\em Combinatorial Enumeration\/}, Wiley, N.Y., 1983.
\bibitem{Kre} G.\ Kreweras, Inversion des polynomes de Bell bidimensionnels et application au denombrement des
    relations binaires connexes. {\em C. R. Acad. Sci.\/} Paris Ser. A-B 268 1969 A577-A579.
\bibitem{Slo1} N.\ J.\ A.\ Sloane, {\em On-Line Encyclopedia of Integer Sequences\/},\\
    {\tt http://www.research.att.com/\~{}njas/sequences/}.
\bibitem{Rio} J.\ Riordan, {\em Combinatorial Identities\/}, John Wiley \& Sons, New York-London-Sydney, 1968.
\bibitem{JK1} V.\ Jovovi\'c and G.\ Kilibarda, On the number of Boolean functions in the Post classes $F^\mu_8$,
    {\em Diskretnaya Matematika\/}, 11 no. 4 (1999), 127--138 (translated in {\em Discrete Mathematics and
    Applications\/}, 9 no. 6, 1999).
\bibitem{Agn} M.\ Aigner, {\em Combinatorial Theory\/}, Springer-Verlag, Berlin Heidelberg New York, 1979.
\bibitem{Slo2} N.\ J.\ A.\ Sloane, {\em Transformations of Integer Sequences\/},\\
    {\tt http://www.research.att.com/\symbol{126}njas/sequences/transforms.html}.
\bibitem{Wlf} H.\ S.\ Wilf, {\em Generatingfunctionology\/}, Academic Press, N.Y., 1994.
\bibitem{Har} F.\ Harary and E.\ M.\ Palmer, {\em Graphical Enumeration\/}, Academic Press, N.Y., 1973.
\bibitem{KJ1} G. Kilibarda and V. Jovovi\'c, On the number of monotone Boolean functions with fixed number of lower
    units (in Russian), {\em Intellektualnye sistemy\/} ({\em Moscow\/}) {\bf 7} (1--4) (2003), 193--217.
\bibitem{KJ2} G. Kilibarda and V. Jovovi\'c, Antichains of Multisets, {\em Journal of Integer Sequences\/},
    Vol.\ 7 (2004), Article 04.1.5.
\end{thebibliography}
\end{document}